\documentclass[12pt]{amsart}
\usepackage{amsmath}
\usepackage{amssymb}
\usepackage{amsmath,amssymb,amsthm,amscd}
\usepackage{geometry}
\usepackage{txfonts}
\usepackage{amsbsy}
\usepackage{graphicx,color}

\numberwithin{equation}{section}

\newtheorem{prop}{Proposition}[section]
\newtheorem{theo}{Theorem}[section]
\newtheorem{lemm}{Lemma}[section]
\newtheorem{coro}{Corollary}[section]

\def\begeq{\begin{equation}}
\def\endeq{\end{equation}}

\begin{document}

\title{Finite time blowup and type II rate for harmonic heat flow from Riemannian manifolds}
\author{Shi-Zhong Du}
\thanks{The author is partially supported by NSFC (12171299), and GDNSF (2019A1515010605)}
  \address{The Department of Mathematics,
            Shantou University, Shantou, 515063, P. R. China.} \email{szdu@stu.edu.cn}

\renewcommand{\subjclassname}{%
  \textup{2010} Mathematics Subject Classification}
\subjclass[2010]{Primary 35A20, Secondary 58E20}
\date{Dec. 2021}
\keywords{Harmonic heat flow, Type II blow up, Finite time blowup.}

\begin{abstract}
   Let $(M,g)$ be a $m$ dimensional Riemannian manifold with metric $g$ and $(N,h)$ be a $n$ dimensional Riemannian sub-manifold of ${\mathbb{R}}^k$ with induced metric $h$. In this paper, we will study the existence of finite time singularity to harmonic heat flow
       $$
        u_t-\triangle_g u=A_u(\nabla u,\nabla u)
       $$
   and their formation patterns.

    After works of Coron-Ghidaglia \cite{CG}, Ding \cite{D} and Chen-Ding \cite{CD}, one knows blow-up solutions under smallness of initial energy for $m\geq3$. Soon later, 2 dimensional blowup solutions were found by Chang-Ding-Ye in \cite{CDY}. The first part of this paper is devoted to construction of new examples of finite time blow-up solutions without smallness conditions for $3\leq m<7$. In fact, when considering rotational symmetric harmonic heat flow from $B_1\subset{\mathbb{R}}^m$ to $S^m\subset{\mathbb{R}}^{m+1}$, we will prove that the maximal solution blows up in finite time if $b>\vartheta_m$, and exists for all time if $0<b<\frac{\pi}{2}$. This result can be regarded as a generalization of results of Chang-Ding-Ye \cite{CDY} and Chang-Ding \cite{CD2} to higher dimensional case, which relies on a completely different argument. The second part of the paper study the rate of blow-up solutions. When $M$ is a bounded domain in ${\mathbb{R}}^2$ and consider Dirichlet boundary condition on $\partial M$, Hamilton \cite{CDY} has obtained that the blowup rate must be faster than $(T-t)^{-1}$. Under a similar setting, it was later improved a little by Topping \cite{T} to $(T-t)^{-1}|\log(T-t)|$. In this paper, we will extend the results to all Riemannian surfaces $M$ and improve the rate of Topping to $(T-t)^{-1}a(|\log(T-t)|)$ for any positive nondecreasing function $a(\tau)$ satisfying
      $$
        \int^\infty_1\frac{d\tau}{a(\tau)}=+\infty,
      $$
   which is comparable to a recent result of Rapha\"{e}l-Schweyer \cite{RS} for rotational symmetric solutions. Turning to the higher dimensional case $3\leq m<7$, we will demonstrate a completely different phenomenon by showing that all rotational symmetric blow-up solutions can not be type II, which is different to the result of $m\geq7$ by Bizo\'{n}-Wasserman \cite{BW}. Finally, we also present result of finite time type I blowup for heat flow from $S^m$ to $S^m\subset{\mathbb{R}}^{m+1}$, when $3\leq m<7$ and degree is no less than $2$.
\end{abstract}

\maketitle\markboth{Shi-Zhong Du}{Finite time blowup and type II rate}

\tableofcontents

\vspace{20pt}

\section{Introduction}

 Let $(M,g)$ be a $m$ dimensional compact Riemannian manifolds with metric $g$, and $(N,h)$ be a $n$ dimensional submanifold (without boundary) of Euclidean space $R^k$ with induced metric $h$. Taking any mapping $u=(u^1,\cdots, u^k)$ from $M$ to $N\subset R^k$, we can define its Dirichlet energy by
    $$
      E(u)=\frac{1}{2}\int_{M}e(u)dV_g,
    $$
 where
              $$
                dV_g=\sqrt{det(g)}dx
              $$
stands for the volume element of $M$ and
    $$
      e(u)\equiv|\nabla u|^2=g^{ij}\frac{\partial u}{\partial x^i}\frac{\partial u}{\partial x^j}
    $$
denotes the density of the energy. Hereafter, we omit the summations for multiple indexes as usually.

  A harmonic map $u: (M,g) \to (N,h)$ is a critical point of the energy $E(\cdot)$ under variation $u_\varepsilon :\ (M,g)\to(N,h)$ for $\varepsilon\in[0,1]$ which satisfies
        $$
          \frac{d}{d\varepsilon}\Bigg|_{\varepsilon=0}u_\varepsilon(x)=\Phi(x),
        $$
  where $\Phi\in T_uN$ is a tangential fields of $N$. Therefore, a mapping $u\in C^2(M,N)$ is harmonic if and only if it satisfies the Euler-Lagrange equation
    \begin{equation}\label{e1.1}
      \tau(u)\equiv \triangle_g u-A_u(\nabla u, \nabla u)=0,
    \end{equation}
  where $\tau(u)$ is the torsion field of $\triangle u$ and $A_u: T_uN\times T_uN\to(T_uN)^\perp$ is the second fundamental form of $N\subset R^k$ at $u$.

  It's natural to consider the heat version of \eqref{e1.1} by
     \begin{equation}\label{e1.2}
       \frac{\partial u}{\partial t}=\tau(u)=\triangle_g u-A_u(\nabla u,\nabla u).
     \end{equation}
  It's well known that for any given smooth initial mapping $u_0\in C^\infty(M,N)$, there exists a local regular solution to \eqref{e1.2} on $\Omega\times[0,\omega)$. It is maximal in the sense
      $$
        \lim_{t\to \omega^-}\max_{x\in M}|\nabla u|^2(x,t)=+\infty
      $$
when $\omega<+\infty$. In this case, we call it to be blowup in finite time.

     A pioneering work of Eells and Sampson \cite{ES} shows that finite time blowup does not occur when target manifold is of nonpositive curvature. For general target manifolds, a natural question arosed: whether finite time singularity can develop? Rather restrictive examples are known until now. A first example was given by Coron-Ghidaglia in \cite{CG}, where they showed that for Cauchy problem \eqref{e1.2} with $m\geq3$, some rotational symmetric solutions from ${\mathbb{R}}^m$ to $S^m$ with small initial energy blows up in finite time. Later, Ding provided more general examples of blow-up solutions in \cite{D} for $3$ dimensional case, assuming only initial mapping belongs to a nontrivial homotopy class and the initial energy is sufficiently small. This result was later generalized by Chen-Ding in \cite{CD} to higher dimensional case. All these results need the smallness of initial energy.

     In case of dimension two, Chang-Ding \cite{CD2} considered rotational symmetric harmonic heat flow from planar disk and proved that the solution must exist for all time when $b<\pi$. Later, Chang-Ding-Ye showed in \cite{CDY} that finite time blowup does occur for $b>\pi$ under a similar setting (see also \cite{DPW} for non-radial symmetric harmonic heat flow from surface which blows up in finite time). The first purpose of this paper is to extend the results of Chang-Ding-Ye \cite{CDY} and Chang-Ding \cite{CD2} to higher dimensional case, by showing that the solution blows up in finite time if  $b>\vartheta_m$, and exists for all time if $0<b<\frac{\pi}{2}$. The main obstacle of this generalization is the difficulty of finding explicitly blowup sub-solutions. Fortunately, using a complete different argument based on intersection comparison, we can prove a desired finite time blowup and long time existence result for $m\geq3$. It's also notable that the latter part of long time existence was a consequence of a result of Jost in \cite{J}, where he showed that the harmonic heat flow in any dimension does not blow up when the image of the flow is contained inside a strictly convex part of the target. A first main purpose of this paper is to show the following finite time blowup and long time existence results.

 \begin{theo}\label{t1.1}
  Suppose that $3\leq m<7$. We have

  (1) If $b>\vartheta_m$, where $\vartheta_m\in(\frac{\pi}{2},\pi)$ is given by Lemma \ref{p5.1}, then all solutions of
      \begin{equation}\label{e1.3}
     \begin{cases}
       \theta_t=\theta_{rr}+\frac{m-1}{r}\theta_r-\frac{m-1}{r^2}\sin \theta\cos \theta, & 0<r<R, t>0\\
       \theta(0,t)=0, \theta(R,t)=b>0, & t>0\\
       \theta(r,0)=\theta_0(r).
     \end{cases}
   \end{equation}
   blow up in finite time.

  (2) If $0<b<\frac{\pi}{2}$, then all solutions of \eqref{e1.3} exist for all time.
\end{theo}

      Due to the presence of blow-up solutions in different settings, it's important to characterize the formation patterns of their singularities. Like that in curvature flow, we can divide the finite singular time $T$ into type I and II as following:

      Suppose that there exists some positive constant $C$, such that
         \begin{equation}\label{e1.4}
            ||e(u(\cdot,t))||_{C^0(M)}\leq\frac{C}{T-t},\ \ \forall t\in[0,T),
         \end{equation}
      we will call the blow-up time to be type I. Otherwise, it is called to be type II. More precisely, letting $T$ be the blowup time and taking any ${p_0}\in M$, $({p_0},T)$ is called type II singular point provided
    \begin{equation}\label{e1.5}
      \limsup_{t\to T^-}(T-t)||e(u(\cdot,t))||_{C^0(D_\delta({p_0}))}=+\infty
    \end{equation}
for any $\delta>0$, where $D_\delta({p_0})\subset M$ denotes the geodesic ball centered at ${p_0}$ and of radius $\delta$.

  Very few result are known except a result by Hamilton \cite{CDY} for planar domain carrying Dirichlet boundary condition. It was mentioned there that every singular point $(x_0,T)$ for harmonic heat flow from a two dimensional bounded domain  must be type II in sense of \eqref{e1.5}. Later, the blowup rate was improved by Topping \cite{T} a little to $(T-t)^{-1}|\log(T-t)|$ under a similar setting. In this paper, we will extend the results to all Riemannian surfaces and improve blowup rate estimate as following.

\begin{theo}\label{t1.2} (Refined type II blowup)
 Letting $(M,g)$ be a compact surface equipping with metric $g$ and $u$ be a maximal solution of \eqref{e1.2} on $M\times[0,\omega)$, we take any positive nondecreasing function $a(\cdot)$ satisfying
      \begin{equation}\label{e1.6}
        \int^\infty_1\frac{ds}{a(s)}=+\infty
     \end{equation}
and set
     $$
      b(t)=a^{-1}(|\log(\omega-t)|).
     $$
Then for any blow-up point $p_0\in M$, there holds
    \begin{equation}\label{e1.7}
      \limsup_{t\to T^-}b(t)(T-t)||e(u(\cdot,t))||_{C^0(D_{\sqrt{b(t)(T-t)}}(p_0))}=+\infty.
    \end{equation}
\end{theo}

   Specially, when $a(s)=s$, the rate matches exactly the one claimed by Topping. As mentioned in \cite{T}, for any $\delta>0$, there is a heat flow from surfaces such that
      $$
        \limsup_{t\to T^-}\lambda(t)||e(u(\cdot,t))||_{C^0\big(D_{\sqrt{\lambda(t)}}(p_0)\big)}<+\infty
      $$
   with
      $$
        \lambda(t)\equiv(T-t)^{1+\delta}.
      $$
So, it's natural to ask whether there exist some counterexamples for any positive nondecreasing function $a(\cdot)$ satisfying
      $$
       \int^\infty_1\frac{ds}{a(s)}<+\infty?
      $$
As shown in a recently paper \cite{RS} by Rapha\"{e}l-Schweyer, this may not be true since the blowup rate \eqref{e1.7} can be improved to a sharp one $a(s)=s^2$ for 2 dimensional rotational symmetric harmonic heat flow. The possibility of the extendable of Rapha\"{e}l-Schweyer's theorem to non-rotational symmetric case would be an interesting problem.

    When considering higher dimensional case, the situation changes dramatically. Actually, we have the following characterization result of blowup rate for $3\leq m<7$.

\begin{theo}\label{t1.3}
  Let $\theta$ be the maximal solution of \eqref{e1.3} on $[0,1]\times[0,\omega)$ with $0<\omega<+\infty$. If $3\leq m<7$ we have
     \begin{equation}\label{e1.8}
         \limsup_{t\to\omega^-}(\omega-t)\sup_{0<r\leq1}\Bigg(\frac{m-1}{r^2}\sin^2\theta+\theta_r^2\Bigg)(r,t)<+\infty.
     \end{equation}
\end{theo}

 As shown by Bizo\'{n}-Wasserman \cite{BW}, for $m\geq7$, all blowups are of type II. Thus the range $m<7$ in Theorem \ref{t1.3} can not be improved in general. Combining our type II Theorem \ref{t1.3} with the above finite time blowup Theorem \ref{t1.1}, new examples of blow-up solutions are derived on higher dimensional case without smallness of initial energy. Furthermore, the blowup rates are determined to be type I, which are new so far.

     At the end of the paper, we will also present the following result of type I finite time blowup for heat flow from $S^m, 3\leq m<7$ to $S^m\subset{\mathbb{R}}^{m+1}$, whose degree is no less than two.

 \begin{theo}\label{t1.4}
  Consider the rotatory-inversion symmetric harmonic heat flow from $S^m$ to $S^m\subset{\mathbb{R}}^{m+1}, 3\leq m<7$. If $u_0$ is a super-harmonic map whose degree is no less than $2$, the solution blows up in finite time with type I rate.\\
\end{theo}

 Our Theorem \ref{t1.4} makes a striking difference from that of \cite{CDZ} for degree is no less than three in dimension two. Contents of this paper are organized as follows. At first, we will give the proof of Theorem \ref{t1.2} in Section 2 and 3 by self-similar variables and a $\varepsilon-$regularity lemma in spirit of Struwe \cite{S1,S2}. Using a Sturm-Liouville type zero comparison principle stated in Section 5, finite time blowup Theorem \ref{t1.1} will be shown in Section 4. Next, we complete the proof of type I rate Theorem \ref{t1.3} on disk after Section 5 and 6. Finally, a proof of Theorem \ref{t1.4} will be presented in Section 7.

\vspace{40pt}

\section{Self-similar variables and decaying local energy}

In this section, we consider a maximal harmonic heat flow of \eqref{e1.2} on compact Riemannian surface $M$ equipping with metric $g$ and impose initial condition
   $$
     u(x,0)=u_0(x)\in N, \ \ \forall x\in M.
   $$
We shall prove two key localized monotonicity formulas of Giga-Kohn type \cite{GK}. (see also \cite{CDZ} for another version) Multiplying \eqref{e1.2} with $u_t$, integrating over $M$ and performing integration by parts, we can derive an energy identity
  \begin{equation}\label{e2.1}
    \frac{d}{dt}E(u)=-\int_\Omega u_t^2dx
  \end{equation}
for a non-increasing energy
   $$
     E(u)\equiv\frac{1}{2}\int_\Omega|\nabla u|^2dx,
   $$
and dissipation of torsion field
  \begin{equation}\label{e2.2}
    \int^\omega_0\int_\Omega u_t^2dxdt\leq E(u_0)<\infty.
  \end{equation}

Throughout this paper, we use
   $$
     D_R(p_0)\equiv\Big\{p\in M\Big| \ dist(p,p_0)\leq R\Big\}, \ \ p_0\in M
   $$
and
   $$
     B_R(x_0)\equiv\Big\{x\in{\mathbb{R}}^m\Big| |x-x_0|\leq R\Big\}
   $$
to denote geodesic ball in $M$ and Euclidean ball in local chart of coordinates respectively. Likewise, we use
   $$
    P_R(p_0,t_0)\equiv\Big\{(p,t)\in M\times{\mathbb{R}}\Big|\ p\in D_R(p_0), \ \ t_0-R^2<t<t_0\Big\}
   $$
for each $(p_0,t_0)\in M\times{\mathbb{R}}$ to denote the usual parabolic cylinder and use
   $$
    Q_R(x_0,t_0)\equiv\Big\{(x,t)\in M\times{\mathbb{R}}\Big|\ x\in B_R(x_0), \ \ t_0-R^2<t<t_0\Big\}
   $$
to denote its counterpart in local coordinates. Since any small neighborhood of $M$ is approximating isometric to a small ball in Euclidean space, without loss of generality, one may assume there exists a positive constant $\beta=\beta_{M,g}>1$ such that
    \begin{equation}\label{e2.3}
       D_{\beta^{-1}R}(p_0)\subset\exp_p\Big\{B_R(x_0)\Big\}\subset D_{\beta R}(p_0)
    \end{equation}
holds for all $p_0\in M$ and $0<R<R_0$, where $x_0$ is the coordinate of $p_0$.

For any point $p_0\in M$, suppose that $M$ can be parameterized by local coordinates $x=(x^1,x^2)\in B_\delta\subset{\mathbb{R}}^2$ in a geodesic ball $D_\delta(a)$ for some $\delta>0$. Without loss of generality, we may assume the coordinate of $p_0$ is $(0,0)$ and
    $$
      g_{ij}(0)=\delta_{ij}, \nabla g_{ij}=0.
    $$
(For example, we can take a normal coordinate here)

Now, taking a cut-off function $\phi=\phi(x)$ supporting in $B_\delta$ as above, we may assume that \eqref{e2.7} and \eqref{e2.8} hold for some positive constant $C$.

 Re-scaling $u$ by self-similar variables
    $$
      x=e^{-\frac{s}{2}}y,\ t=\omega-e^{-s}
    $$
 and setting
    $$
     u(x,t)=w\Bigg(\frac{x}{\sqrt{\omega-t}},\ -\log(\omega-t)\Bigg),
    $$
 the function $w=w_{(p_0,\omega)}(y,s)$ satisfies
    \begin{equation}\label{e2.4}
      w_s-\triangle_{\widetilde{g}}w+\frac{1}{2}y_i\nabla_iw
      =A_w(\nabla_{\widetilde{g}}w,\nabla_{\widetilde{g}}w),
    \end{equation}
 where
    $$
     \widetilde{g}_{ij}(y,s)\equiv g_{ij}(e^{-\frac{s}{2}}y+x_0).
    $$

 Noting that the Laplace-Beltrami operator is given by
   $$
     \triangle_{\widetilde{g}}=\frac{1}{\sqrt{\widetilde{g}}}\frac{\partial}{\partial x^i}\Big(\sqrt{\widetilde{g}}\ \widetilde{g}^{ij}\frac{\partial}{\partial x^j}\Big)
   $$
with $\widetilde{g}=det(\widetilde{g}_{ij})$, we can rewrite the principle part of \eqref{e4.2} into self-adjoint form
   \begin{equation}\label{e2.5}
      \rho w_s-\nabla_i(\rho\widetilde{g}^{ij}\nabla_jw)-\frac{\rho}{2}(\widetilde{g}^{ij}-\delta_{ij})y_i\nabla_jw=\rho A_w(\nabla w,\nabla w),
   \end{equation}
where $\nabla_i=\nabla_{\frac{\partial}{\partial y^i}}$ denotes the Levi-Civita connection on $M$ and $\rho(y)=e^{-\frac{|y|^2}{4}}$. For any $\delta>0$, we assign a cut-off function
   \begin{equation}\label{e2.6}
     \phi(x)=\begin{cases}
        1, & x\in B_{\delta/2},\\
        0, & x\not\in B_\delta
     \end{cases}
   \end{equation}
to the given point $p_0\in M$, which satisfies that
   \begin{equation}\label{e2.7}
    0\leq\phi(x)\leq1, \ \ |\nabla\phi|(x)\leq C_*\delta^{-1},\\ |\nabla^2\phi|(x)\leq C_*\delta^{-2},\ \ \forall x\in{\mathbb{R}}^2
   \end{equation}
with universal constant $C_*$. Expressing in terms of self-similar variables, one gets an expanding cut-off function
   $$
     \varphi(y,s)=\phi(e^{-\frac{s}{2}}y)
   $$
satisfying that
   \begin{eqnarray}\label{e2.8}\nonumber
    |\nabla\varphi(y,s)|&\leq& Ce^{-\frac{s}{2}}\chi_{B_{\delta e^{s/2}}\setminus B_{\frac{\delta}{2}e^{s/2}}},\ \ \ |\nabla^2\varphi(y,s)|\leq Ce^{-s}\chi_{B_{\delta e^{s/2}}\setminus B_{\frac{\delta}{2}e^{s/2}}},\\
    &&\ \ \ \ |\varphi_s(y,s)|\leq C\chi_{B_{\delta e^{s/2}}\setminus B_{\frac{\delta}{2}e^{s/2}}}
   \end{eqnarray}
for all $y\in{\mathbb{R}}^2$. Hereafter, we denote
  $$
   |\nabla\varphi|_{\widetilde{g}}^2\equiv\widetilde{g}^{ij}\nabla_i\varphi\nabla_j\varphi,\ \ \ |\nabla\varphi|^2\equiv\Sigma_{i=1}^m\Bigg(\frac{\partial\varphi}{\partial x^i}\Bigg)^2
  $$
and use the fact that
   $$
     \frac{1}{C_{M}}|\nabla\varphi|^2\leq|\nabla\varphi|_{\widetilde{g}}^2\leq C_{M}|\nabla\varphi|^2
   $$
holds for some positive constant $C_{M}$ depending only on $M$.

 Multiplying \eqref{e4.3} by $w_s\varphi^2$, integrating over $M$ and performing integration by parts, one gets
   \begin{eqnarray}\label{e2.9}\nonumber
     &\int_{M}\rho w_s^2\varphi^2dV_{\widetilde{g}}+\frac{1}{2}\frac{d}{ds}\int_{M}\rho|\nabla w|_{\widetilde{g}}^2\varphi^2dV_{\widetilde{g}}=\frac{1}{2}\int_{M}\rho\frac{\partial}{\partial s}\widetilde{g}^{ij}\nabla_iw\nabla_jw\varphi^2dV_{\widetilde{g}}+\int_{M}\rho|\nabla w|_{\widetilde{g}}^2\varphi\varphi_sdV_{\widetilde{g}}&\\ \nonumber
     &-2\int_{M}\rho\widetilde{g}^{ij}\nabla_iww_s\varphi\nabla_j\varphi dV_{\widetilde{g}}
     +\frac{1}{2}\int_{M}\rho(\widetilde{g}^{ij}-\delta_{ij})y_i\nabla_jww_s\varphi^2dV_{\widetilde{g}}&\\ \nonumber
     &\leq\varepsilon\int_{M}\rho w_s^2\varphi^2dV_{\widetilde{g}}+C_\varepsilon\Bigg(e^{-\frac{s}{2}}\int_{M}\rho|y||\nabla w|_{\widetilde{g}}^2\varphi^2dV_{\widetilde{g}}+e^{-s}\int_{M}\rho|y|^4|\nabla w|_{\widetilde{g}}^2\varphi^2dV_{\widetilde{g}}\Bigg)&\\
     &+C_\varepsilon\int_{M}\rho|\nabla w|_{\widetilde{g}}^2(\varphi|\varphi_s|+|\nabla\varphi|^2)dV_{\widetilde{g}},&
   \end{eqnarray}
for any small $\varepsilon>0$ and large $C_\varepsilon$, where
    \begin{eqnarray*}
      \Bigg|\frac{\partial}{\partial s}\widetilde{g}^{ij}\Bigg|&=&\frac{1}{2}e^{-\frac{s}{2}}|y\cdot\nabla_xg_{lk}\widetilde{g}^{il}\widetilde{g}^{jk}|\\
      &\leq& Ce^{-\frac{s}{2}}|y|
    \end{eqnarray*}
and
    \begin{eqnarray*}
      |\widetilde{g}^{ij}-\delta_{ij}|\leq Ce^{-s/2}|y|
    \end{eqnarray*}
have been used. Another hand, for any given $\kappa>1$, there exists a large number $C_\kappa$, such that
   \begin{eqnarray}\label{e2.10}\nonumber
     e^{-\frac{s}{2}}\rho|y|&\leq& C_\kappa\Big(e^{-\frac{s}{4}}\rho+e^{-\kappa s}\Big)\\
     e^{-s}\rho(|y|^2+|y|^4)&\leq& C_\kappa\Big(e^{-\frac{s}{4}}\rho+e^{-\kappa s}\Big).
   \end{eqnarray}
So, it yields from \eqref{e2.9} \eqref{e2.10} and \eqref{e2.9} that
  \begin{equation}\label{e2.11}
     \frac{d}{ds}{\mathcal{E}}(w)\leq-(1-\varepsilon)\int_{M}\rho w_s^2\varphi^2dV_{\widetilde{g}}+C_{\kappa,\varepsilon} e^{-\frac{s}{4}}\int_{M}\rho|\nabla w|_{\widetilde{g}}^2\varphi^2dV_{\widetilde{g}}+C_{\kappa,\varepsilon} e^{-\kappa s}E(u_0),
  \end{equation}
for
   $$
    {\mathcal{E}}(w)=\frac{1}{2}\int_{M}\rho|\nabla w|_{\widetilde{g}}^2\varphi^2dV_{\widetilde{g}}.
   $$

Integrating over time, we derive the first key monotonicity formula:

\begin{prop}\label{p2.1}
  Assume that $(M,g)$ is a compact Riemannian surface equipped with metric $g$. Let $p_0\in M$ and $w=w^{p_0}(y,s)$ be the rescaled solution of \eqref{e2.5}. Then for any $\kappa>1$, there exist  two large constants $C_{M,\kappa}$ and $s_{M,\kappa}$ depending only on $M$ and $\kappa$, such that
    \begin{equation}\label{e2.12}
      {\mathcal{E}}(w(s))\leq(1+\kappa^{-1}){\mathcal{E}}(w(s'))+C_{M,\kappa} e^{-\kappa s'}E(u_0)\ \ \forall s>s'\geq s_{M,\kappa}
    \end{equation}
  and
    \begin{equation}\label{e2.13}
      \int^\infty_{-\log\omega}\int_{M}\rho w_s^2\varphi^2dV_{\widetilde{g}}dt\leq C_{M,\kappa}  E(u_0)<\infty.
    \end{equation}
\end{prop}

\noindent\textbf{Remark.} Original version of this monotonicity formula \eqref{e2.12} for harmonic heat flow was firstly found by Struwe in  \cite{S1}. It was later developed in \cite{S2} and subsequential literatures. We present here a localized version for the purpose of proving our type II rate result.\\

\vspace{10pt}

  To explore the structure of singular point $(p_0,\omega)$ further, we need to introduce a second local energy by
   $$
    \widetilde{{\mathcal{E}}}(w)\equiv\frac{1}{2}\int_{M}\rho|y|^2|\nabla w|_{\widetilde{g}}^2\varphi^2dV_{\widetilde{g}}
   $$
and prove the following key decaying estimate.

\begin{prop}\label{p2.2}
   Assume that $(M,g)$ is a compact Riemannian surface equipped with metric $g$. Let $p_0\in M$ and $w=w^{p_0}(y,s)$ be the rescaled solution of \eqref{e2.5}. Then for any $\lambda>8, \kappa>1$, there exist two large constants $C'_{\lambda,\kappa}$ and $s'_{\lambda,\kappa}$ depending only on $\lambda$ and $\kappa$, such that
     \begin{eqnarray}\label{e2.14}\nonumber
       \widetilde{{\mathcal{E}}}(w(s))+\lambda{\mathcal{E}}(w(s))&\leq&(1+\kappa^{-1})\Big(\widetilde{{\mathcal{E}}}(w(s'))+\lambda{\mathcal{E}}(w(s'))\Big)\\
        &&+C_{\lambda,\kappa} e^{-\kappa s}E(u_0),\ \ \forall s>s'\geq s'_{\lambda,\kappa}
     \end{eqnarray}
   and
     \begin{eqnarray}\label{e2.15}\nonumber
      && \int^\infty_{-\log\omega}\int_{M}\Bigg(\rho(|y|^2+1)w_s^2\varphi^2+\rho|y|^2|\nabla w|_{\widetilde{g}}^2\varphi^2\Bigg)dV_{\widetilde{g}}ds\\
      && \ \ \ \ \ \ \ \ \ \ \ \ \  \ \ \ \ \ \ \ \ \ \ \ \ \  \ \ \ \ \ \ \ \ \ \ \ \ \  \leq C_{\lambda,\kappa} E(u_0).
     \end{eqnarray}
 \end{prop}

 \noindent\textbf{Proof.} Multiplying \eqref{e4.3} by $|y|^2w_s\varphi^2$, integrating over $M$ and performing integration by parts, one gets
  \begin{eqnarray}\label{e2.16}\nonumber
    \frac{d}{ds}\widetilde{{\mathcal{E}}}(w)&\leq&-(1-\varepsilon)\int_{M}\rho|y|^2w_s^2\varphi^2dV_g-2\int_{M}\rho\Big<y,\nabla w\Big>_{\widetilde{g}}w_s\varphi^2dV_{\widetilde{g}}\\
    &&+C_{M,\varepsilon,\kappa}e^{-\frac{s}{4}}\widetilde{{\mathcal{E}}}(w)+C_{M,\varepsilon,\kappa}e^{-\kappa s}E(u_0)
  \end{eqnarray}
for some positive constant $C_{M,\varepsilon,\kappa}$ depending only on $M$ and $\varepsilon>0, \kappa>1$, where
  \begin{eqnarray*}
      \Big<u,v\Big>_{\widetilde{g}}&\equiv&\widetilde{g}^{ij}u_iv_j,\ \ \ |u|^2_{\widetilde{g}}\equiv\Big<u,u\Big>_{\widetilde{g}}
  \end{eqnarray*}
and \eqref{e2.10} has been used again. To control the second term on R.H.S. of \eqref{e2.16}, we need to drive the Pohozaev identity as before: multiplying \eqref{e1.2} by $\Big<y,\nabla w\Big>_{\widetilde{g}}\varphi^2$ and performing integration by parts, we obtain that for any $\varepsilon>0, \kappa>1$, there exist two positive constants $s_{M,\varepsilon, \kappa}, C_{M,\kappa}$, such that
  \begin{eqnarray}\label{e2.17}\nonumber
    &\int_{M}\rho\Big<y,\nabla w\Big>_{\widetilde{g}}w_s\varphi^2dV_{\widetilde{g}}=-\int_{M}\rho\widetilde{g}^{ij}\nabla_jw\Bigg(\nabla_i\widetilde{g}^{kl}y_k\nabla_lw\varphi^2+\widetilde{g}^{il}\nabla_lw\varphi^2+\widetilde{g}^{kl}y_k\frac{\partial^2w}{\partial y^i\partial y^l}\varphi^2\Bigg)dV_{\widetilde{g}}&\\ \nonumber
    &-2\int_{M}\rho\Big<\nabla w,\nabla\varphi\Big>_{\widetilde{g}}\Big<y,\nabla w\Big>_{\widetilde{g}}\varphi dV_{\widetilde{g}}+\frac{1}{2}\int_{\Omega_s}\rho(\widetilde{g}^{ij}-\delta_{ij})y_i\nabla_jw\Big<y,\nabla w\Big>_{\widetilde{g}}\varphi^2dV_{\widetilde{g}}&\\ \nonumber
    &\leq-\frac{1}{4}\int_{M}\rho|y|^2|\nabla w|_{\widetilde{g}}^2\varphi^2dV_{\widetilde{g}}+C_{M,\kappa}e^{-\frac{s}{2}}{\mathcal{E}}(w(s))&\\
    &+C_{M,\kappa} e^{-\frac{s}{2}}\int_{M}\rho|y|^2|\nabla w|_{\widetilde{g}}^2\varphi^2+C_{M,\kappa} e^{-\kappa s}E(u_0)&\\ \nonumber
    &\leq-\Big(\frac{1}{4}-\varepsilon\Big)\int_{M}\rho|y|^2|\nabla w|_{\widetilde{g}}^2\varphi^2dV_{\widetilde{g}}+C_{M,\kappa}e^{-\frac{s}{2}}{\mathcal{E}}(w(s))+C_{M,\kappa} e^{-\kappa s}E(u_0)&
  \end{eqnarray}
for all $s>s_{M,\varepsilon,\kappa}$. Combining \eqref{e2.16} \eqref{e2.17} with \eqref{e2.11}, the conclusion follows from the positivity of
   \begin{eqnarray*}
    \lambda\int_M\rho w_s^2\varphi^2dV_{\widetilde{g}}+(2+\mu)\int_M\rho\langle y,\nabla w\rangle_{\widetilde{g}}w_s\varphi^2dV_{\widetilde{g}}+\frac{\mu}{4}\int_M\rho|y|^2|\nabla w|^2_{\widetilde{g}}\varphi^2dV_{\widetilde{g}}
   \end{eqnarray*}
for any $\lambda>8$ and some positive $\mu=\mu_\lambda$. $\Box$

\vspace{40pt}

\section{Nondegeneracy of blowup and refined type II rates}

In this section ,we will show the following type II result.

\begin{theo}\label{t3.1}
  Assume that $(M,g)$ is a compact Riemannian surface equipped with a metric $g$ and $u$ is a maximal solution to \eqref{e1.2} on $M\times[0,\omega)$. We have the blowup set is discrete and all blowup points $({p_0},\omega)$ must be refined type II, say,
     \begin{equation}\label{e3.1}
       \limsup_{t\to\omega^-}\sqrt{b(t)(\omega-t)}\sup_{dist(p,p_0)\leq \sqrt{b(t)(\omega-t)}}|\nabla u|(p,t)=+\infty,
     \end{equation}
  holds for
      $$
        b(t)=a^{-1}(|\log(\omega-t)|)
      $$
  and a positive nondecreasing function $a(s)$ satisfying \eqref{e1.6}.
\end{theo}

 In order proving the main theorem \ref{t3.1}, we need several crucial lemmas. The first one is the following Harnack type inequality for linear partial differential inequality.

\begin{lemm}\label{l3.1}
   Let $(M,g)$ be a $m$ dimensional compact manifold with metric $g$ and $0<R\leq 1$. If $v$ be a nonnegative classical solution to
    \begin{equation}\label{e3.2}
     v_t-\triangle_gv\leq0, \ \ \mbox{ in }\overline{P_R({p_0},0)}\subset M,
    \end{equation}
 then
     \begin{equation}\label{e3.3}
       v({p_0},0)\leq C_M\fint_{P_R({p_0},0)}vdV_gdt
     \end{equation}
  holds for some positive constant $C_M$ depending only on $M$.
\end{lemm}

Before proving Lemma \ref{l3.1}, let's first recall a Sobolev inequality \cite{GT} in the following form:

\begin{lemm}\label{l3.2}
  Let $M$ be a $m$-dimensional compact manifold with metric $g$ and $\Omega\subset M$ be any relative open subset of $M$ with smooth boundary $\partial\Omega$. There exists a positive constant $C_g$ depending only on the least and largest eigenvalues of $g$, such that
     \begin{equation}\label{e3.4}
       \int_\Omega |u|^{\frac{2m+4}{m}}dV_g\leq C_g\Bigg(\int_\Omega u^2dV_g\Bigg)^{\frac{2}{m}}\int_\Omega|\nabla u|_g^2dV_g
     \end{equation}
   holds for any $u\in W^{1,2}_0(\Omega)$.\\
\end{lemm}

\noindent\textbf{Proof of Lemma \ref{l3.1}.} For any $q>1$, direct computation shows that
  \begin{eqnarray}\nonumber\label{e3.5}
    (\partial_t-\triangle_g)v^q&=&qv^{q-1}(v_t-\triangle_gv)-q(q-1)v^{q-2}|\nabla v|^2\\
      &\leq&-q(q-1)v^{q-2}|\nabla v|^2.
  \end{eqnarray}
For any $r\in(0,R]$, let's take a cut-off function $\xi\in C^\infty(Q_r)$ satisfying that
   $$
    \xi(x,t)=\begin{cases}
      1, & (x,t)\in B_{r/2}\times[-r^2/4,0],\\
      0, & (x,t)\not\in B_r\times[-r^2,0]
    \end{cases}
   $$
and
   \begin{equation}\label{e3.6}
     0\leq\xi\leq 1, \ \ r|\nabla\xi|_g+r^2(|\xi_t|+|\triangle_g\xi|)\leq C_M
   \end{equation}
for some positive constant $C_M$ depending only on $M$. Multiplying \eqref{e3.5} by $\xi^2$ and integrating over space-time, we obtain that
   \begin{eqnarray}\nonumber\label{e3.7}
    && \sup_{t\in[-r^2,0]}\int_{B_r}v^q\xi^2dV_g+\frac{2(q-1)}{q}\int^0_{-r^2}\int_{B_r}|\nabla(v^{\frac{q}{2}}\xi)|_g^2dV_gdt\\
    && \ \ \ \ \ \ \ \ \ \ \ \ \ \ \leq 4\iint_{Q_r}v^q\Big(\xi|\xi_t|+\varphi|\triangle_g\xi|+|\nabla\xi|_g^2\Big)dV_gdt\leq 4C_{M}r^{-2}\iint_{Q_r}v^qdV_gdt,
   \end{eqnarray}
where
   \begin{eqnarray*}
     &\int_{Q_r}|\nabla(v^{\frac{q}{2}}\xi)|_g^2dV_g=\int_{Q_r}\Big|\frac{q}{2}v^{\frac{q}{2}-1}\nabla v\xi+v^{\frac{q}{2}}\nabla\xi\Big|_g^2dV_g\leq\frac{q^2}{2}\int_{Q_r} v^{q-2}|\nabla v|_g^2\xi^2dV_g+2\int_{Q_r} v^q|\nabla\xi|_g^2dV_g&\\
     &\Rightarrow \int_{Q_r} v^{q-2}|\nabla v|_g^2\xi^2dV_g\geq\frac{2}{q^2}\int_{Q_r}|\nabla(v^{\frac{q}{2}}\xi)|_g^2dV_g-\frac{4}{q^2}\int_{Q_r} v^2|\nabla\xi|_g^2dV_g&
   \end{eqnarray*}
has been used. Applying Lemma \ref{l3.2} to $u=v^{\frac{q}{2}}\xi$, one can verify that
   \begin{equation}\label{e3.8}
     \Bigg(\fint_{Q_{r/2}}v^{\frac{q(m+2)}{m}}dV_gdt\Bigg)^{\frac{m}{q(m+2)}}\leq (C_M)^{\frac{m}{q(m+2)}}\Bigg(\fint_{Q_r}v^qdV_gdt\Bigg)^{\frac{1}{q}}
   \end{equation}
from \eqref{e3.7}. Iterating $r=R/\beta\cdot 2^{-k}$ and $q=q'\cdot\Big(\frac{m+2}{m}\Big)^k, q'>1$ for $k\geq0, k\in{\mathbb{Z}}$, we obtain that
   \begin{eqnarray}\nonumber\label{e3.9}
     \log a_{k+1}&\leq&\log a_k+\Bigg(\frac{m}{m+2}\Bigg)^{k+1}\Bigg(\frac{1}{q'}\log C_M\Bigg)\\
      &\leq& \log a_0+\Bigg(1-\Bigg(\frac{m}{m+2}\Bigg)^{k+2}\Bigg)\frac{m+2}{2q'}\log C_M,
   \end{eqnarray}
where
   $$
    a_k\equiv\Bigg(\fint_{Q_{R/\beta\cdot 2^{-k}}}v^{q'\cdot\big(\frac{m+2}{m}\big)^k}dV_gdt\Bigg)^{\frac{1}{q'}\big(\frac{m}{m+2}\big)^k}.
   $$
Hence, we get
   \begin{equation}\label{e3.10}
     v(0,0)\leq C_M\Bigg(\fint_{Q_{R/\beta}}v^{q'}dV_gdt\Bigg)^{\frac{1}{q'}}\leq C_M\Bigg(\fint_{P_{R}}v^{q'}dV_gdt\Bigg)^{\frac{1}{q}}, \ \ \forall q>1
   \end{equation}
by sending $k\to+\infty$. Passing to the limit $q'\to 1^+$ in \eqref{e3.10}, we complete the proof of Lemma \ref{l3.1}. $\Box$\\

\begin{lemm}\label{l3.3}
  Let $u$ be a maximal solution to \eqref{e1.2} on $M\times[0,\omega)$. There exist positive constants $\varepsilon_0, R_0, \gamma$ depending only on $M, N$, such that for any fixed point ${p_0}\in M$, if
     $$
      r^{-m}\iint_{P_r(\overline{z})}|\nabla u|_g^2dV_gdt<\varepsilon_0
     $$
  holds for all cylinders $P_r(\overline{z}), \overline{z}=(\overline{{p_0}},\overline{t})$ contained inside the cylinder $P_R({p_0},\omega)$ with $0<R<R_0$, then $({p_0},\omega)$ is not a singular point and there holds
     $$
      \sup_{P_{\gamma R}({p_0},\omega)}|\nabla u|_g\leq 4\beta R^{-1}.
     $$
\end{lemm}

\noindent\textbf{Proof.} Without loss of generality, we may assume that $u(x,t)$ is smooth up to singular time $t=\omega$. Otherwise, we can shift $\omega$ to $\omega-\frac{1}{i}$ and then let $i\to+\infty$. By \eqref{e2.3}, we need only to prove the conclusion under the assumption
   $$
    r^{-m}\iint_{Q_r(\overline{z})}|\nabla u|_g^2dxdt<\varepsilon_0
   $$
for all $Q_r(\overline{z}), \overline{z}=(\overline{x},\overline{t})$ contained inside the cylinder $Q_{R/\beta}(0,\omega)$. Consider
    $$
      K=\sup_{0<r<R/\beta}\Bigg[(R/\beta-r)\sup_{\overline{Q_r(0,\omega)}}|\nabla u|_g\Bigg]
    $$
 and let $r_0\in[0,R), z^*\in\overline{Q_{r_0}(0,\omega)}$ satisfying
    $$
     K=(R/\beta-r_0)|\nabla u|_g(z^*).
    $$
Setting $r_1=\frac{R/\beta-r_0}{2}$, we have
    $$
     Q_{r_1}(z^*)\subset Q_{R/\beta}(0,\omega)
    $$
 and
    $$
     r_1\sup_{Q_{r_1}(z^*)}|\nabla u|_g\leq K.
    $$
 As a result,
   $$
    \sup_{Q_{r_1}(z^*)}|\nabla u|_g\leq\Bigg(\frac{R/\beta-r_0}{r_1}\Bigg)|\nabla u|_g(z^*)=2|\nabla u|_g(z^*).
   $$
 Re-scaling $u$ by
   $$
    \widetilde{u}(y',s')=u(x^*+\mu^{-1}y', t^*+\mu^{-2}s'),\ \ \ \mu\equiv|\nabla u|_g(z^*),
   $$
 one gets
    \begin{equation}\label{e3.11}
      \begin{cases}
        \widetilde{u}_{s'}=\triangle_{\widetilde{g}}\widetilde{u}+A_{\widetilde{u}}(\nabla\widetilde{u},\nabla\widetilde{u}),\\
        |\nabla\widetilde{u}|_{\widetilde{g}}\leq 2, \ |\nabla\widetilde{u}|_{\widetilde{g}}(0,0)=1
      \end{cases}
    \end{equation}
 for $(y',s')\in Q_{\mu r_1}(0,0)$, where
   $$
    \widetilde{g}(y')\equiv g(x^*+\mu^{-1}y').
   $$
  Consequently, it's inferred from Bochner's identity that
   \begin{eqnarray}\label{e3.12}\nonumber
     \frac{\partial}{\partial s'}|\nabla\widetilde{u}|_{\widetilde{g}}^2&=&\triangle_{\widetilde{g}}|\nabla\widetilde{u}|_{\widetilde{g}}^2-2|\nabla^2\widetilde{u}|_{\widetilde{g}}^2
     +2Ric(\nabla\widetilde{u},\nabla\widetilde{u})+2\nabla\widetilde{u}\cdot\nabla\Big(A_{\widetilde{u}}(\nabla\widetilde{u},\nabla\widetilde{u})\Big)\\
     &\leq&\triangle_{\widetilde{g}}|\nabla\widetilde{u}|_{\widetilde{g}}^2+C_N(|\nabla\widetilde{u}|_{\widetilde{g}}^2+|\nabla\widetilde{u}|_{\widetilde{g}}^4)\leq\triangle|\nabla\widetilde{u}|_{\widetilde{g}}^2+4C_N|\nabla\widetilde{u}|_{\widetilde{g}}^2
   \end{eqnarray}
holds in $P_{\mu r_1/\beta}(p_0,0)$, where $C_N$ is a constant depending on $N$. We claim that $K\leq2\beta$. For, if $K>2\beta$, then
   $$
    \mu r_1/\beta\geq1
   $$
and \eqref{e3.11} \eqref{e3.12} hold in $P_1({p_0},0)$. Thus, we have
  \begin{equation}\label{e3.13}
    \iint_{P_1({p_0},0)}|\nabla\widetilde{u}|_{\widetilde{g}}^4\leq 4\iint_{Q_\beta(0,0)}|\nabla\widetilde{u}|_{\widetilde{g}}^2\leq4\mu^{m}\iint_{P_{\beta^2/\mu}(z^*)}|\nabla u|_g^2<4\beta^2\varepsilon_0.
  \end{equation}
Regarding \eqref{e3.12} as a linear parabolic inequality and applying Lemma \ref{l3.1} to
  $$
    v(y',s')=e^{-4C_Ns'}|\nabla\widetilde{u}|_{\widetilde{g}}^2
  $$
with $q=2$, one gets
  \begin{eqnarray*}
    1&=&|\nabla\widetilde{u}|_{\widetilde{g}}^2({p_0},0)=v({p_0},0)\leq C_m\Bigg(\fint_{P_1({p_0},0)}v^2\Bigg)^{\frac{1}{2}}\\
      &\leq& C_{m,N}\Bigg(\fint_{P_1({p_0},0)}|\nabla\widetilde{u}|_{\widetilde{g}}^4\Bigg)^{\frac{1}{2}}\leq C_{m,N}\varepsilon_0^{\frac{1}{2}}.
  \end{eqnarray*}
Contradiction holds as long as $\varepsilon_0$ is chosen small. So $K\leq2\beta$ and hence
  $$
   \sup_{P_{R/(2\beta)}({p_0},\omega)}|\nabla u|_g\leq 4\beta^2R^{-1}.
  $$
$\Box$\\

Now, we can prove a central $\varepsilon-$regularity type result by utilizing the monotonicity formula Proposition \ref{p2.1}:

\begin{theo}\label{t3.2}
  Let $(M,g)$ be a compact surface equipped with metric $g$ and $u$ be a maximal solution to \eqref{e1.2} on $M\times[0,\omega)$. For any $\delta>0$ and the assigned cut-off function $\varphi$ to ${p_0}$, we introduce the first local energy ${\mathcal{E}}(w)$ of the rescaled solution $w=w_{({p_0},\omega)}(y,s)$ as above. Then there exist two positive constants $\varepsilon_1, s_2\geq s_1$ depending only on $M, N$ and $\delta$, such that
     $$
      {\mathcal{E}}(w(s_*))<\varepsilon_1\ \ \mbox{ for some } s_*>s_2
     $$
  implies $({p_0},\omega)$ is not a singular point.
\end{theo}

\noindent\textbf{Proof.} We will prove that for some $\sigma_0>0$ (depending on $u$ and $s_*$), there holds
  \begin{equation}\label{e3.14}
    \sup_{P_{\frac{\gamma\sigma_0}{\beta}}({p_0},\omega)}|\nabla u|_g\leq 4\beta^2\sigma_0^{-1}.
  \end{equation}
In fact, by smoothness of $u$, before blowup time $\omega$, $\exists0<\sigma_0<\min\Big\{\frac{1}{2}e^{-s_*/2},\delta\Big\}$, such that for any
   $$
    (a',\omega')\in\overline{B_{\sqrt{\frac{m}{2\pi e}}\sigma_0}({p_0})\times(\omega-\sigma_0^2,\omega+\sigma_0^2)},
   $$
we have
  \begin{equation}\label{e3.15}
    {\mathcal{E}}(w_{({p_0}',\omega')}(s'_*))<2\varepsilon_1,
  \end{equation}
where
   $$
     s'_*\equiv-\log(\omega'-\omega+e^{-s_*})\in\Bigg(s_*-\log\frac{5}{4}, \ s_*+\log\frac{4}{3}\Bigg)
   $$
Thus, it follows from monotonicity formula \eqref{e4.7} that
   $$
     {\mathcal{E}}(w_{({p_0}',\omega')}(s))\leq \varepsilon, \ \ \forall s\geq s'_*,
   $$
where $\varepsilon$ is small as long as $\varepsilon_1$ is small. Integrating over time, we have
  \begin{equation}\label{e3.16}
    \int^{s+1}_s\int_{M}\rho|\nabla w_{({p_0}',\omega')}|_{\widetilde{g}}^2\varphi^2dV_gd\tau=2\int^{s+1}_s{\mathcal{E}}(w_{({p_0}',\omega')}(\tau))d\tau<4\varepsilon_1.
  \end{equation}
Writting \eqref{e3.16} back to $u(x,t)$, there holds
   \begin{equation}\label{e3.17}
     r^{-m}\int^{\omega'-\frac{r^2}{e}}_{\omega'-r^2}\int_{D_{\sqrt{\frac{m}{2\pi e}}r}({p_0}')}|\nabla u|_g^2dxdt\leq\varepsilon'\equiv4e^{\frac{m}{8\pi e}}\varepsilon_1,
   \end{equation}
for $r=e^{-s/2}$. Changing of parameters
   $$
    \overline{x}={p_0}', \ \overline{t}=\omega'-\frac{r^2}{e},\ \overline{r}=\sqrt{\frac{e-1}{e}}r
   $$
in \eqref{e3.17}, one gets
   \begin{equation}\label{e3.18}
     (r/\beta)^{-m}\iint_{P_{r/\beta}(\overline{z})}|\nabla u|_g^2dV_gdt\leq Cr^{-m}\iint_{Q_r(\overline{z})}|\nabla u|_g^2dxdt\leq C\varepsilon'\leq\varepsilon_0
   \end{equation}
for any $P_{r/\beta}(\overline{z})$ contained inside $P_{\sigma_0/\beta}({p_0},\omega)$, provided $\varepsilon_1$ is chosen small. Now, the conclusion follows from Lemma \ref{l3.3}. $\Box$\\

 Under below, we take a positive monotone non-increasing function $\alpha(s)$, satisfying either one of the following hypothesises:

 (H1) $\lim_{s\to+\infty}\alpha(s)=0$ and $\alpha(s)\geq a^{-\frac{1}{2}}(s), \forall s>1$, where $a(s)$ is a positive non-decreasing function satisfying
    $$
     \int^\infty_1\frac{ds}{a(s)}=\infty,
    $$
 or

 (H2) $\lim_{s\to+\infty}\alpha(s)=\alpha_0>0$.\\

 By dividing the domain into nearby or faraway ones from origin, we can estimate the local energy as following:

 \begin{lemm}\label{l3.4}
     Assume that $(M,g)$ is a compact surface equipped with metric $g$ and ${p_0}\in M$. Letting $w=w_{({p_0},\omega)}(y,s)$ be a re-scaled solution to \eqref{e4.3}, we have the following estimates on local energy ${\mathcal{E}}(w)$:

 Case1: If $\alpha(s)$ satisfies (H1), then there holds
   \begin{eqnarray}\label{e3.19}\nonumber
     \int^{s+1}_s{\mathcal{E}}(w(\tau))d\tau&\leq& a(s+1)\int^{s+1}_s\int_{M}\rho|y|^2|\nabla w|_{\widetilde{g}}^2\varphi^2dV_{\widetilde{g}}d\tau\\
      &&+\frac{1}{2}\int^{s+1}_s\int_{|y|\leq\alpha(\tau)}|\nabla w|_{\widetilde{g}}^2dV_{\widetilde{g}}d\tau.
   \end{eqnarray}

 Case2: If $\alpha(s)$ satisfies (H2), then there holds
    \begin{eqnarray}\label{e3.20}\nonumber
      \int^{s+1}_s{\mathcal{E}}(w(\tau))d\tau&\leq&\frac{1}{2}\alpha_0^{-2}\Bigg(\int^{s+1}_s\int_{\Omega_\tau}\rho|y|^2|\nabla w|_{\widetilde{g}}^2\varphi^2dV_{\widetilde{g}}d\tau\Bigg)\\
      &&+\frac{1}{2}\int^{s+1}_s\int_{|y|\leq\alpha(\tau)}|\nabla w|_{\widetilde{g}}^2dV_{\widetilde{g}}d\tau.
   \end{eqnarray}
\end{lemm}

\noindent\textbf{Proof.} The lemma is a easy consequence of the following inequality
   $$
    \int^{s+1}_s\int_{|y|>\alpha(\tau)}\rho|\nabla w|_{\widetilde{g}}^2\varphi^2dyd\tau\leq\alpha^{-2}(s+1)\int^{s+1}_s\int_{|y|>\alpha(\tau)}\rho|y|^2|\nabla w|_{\widetilde{g}}^2\varphi^2dV_{\widetilde{g}}d\tau.
   $$
$\Box$\\

As a corollary of the decaying estimate Proposition \ref{p2.2} and the $\varepsilon-$regularity Theorem \ref{t3.2}, together with Lemma \ref{l3.4}, we have the following property:

\begin{prop}\label{p3.1}
  Let $(M,g)$ be a compact surface equipped with metric $g$ and $u$ be a maximal solution of \eqref{e1.2} on $M\times[0,\omega)$. Then for any ${p_0}\in M$, there exist two positive constants $\varepsilon_2, s_3\geq s_2$ depending only on $M$ and $N$, such that $({p_0},\omega)$ is not a singular point, provided

Case1: $\alpha(s)$ satisfies (H1) and
  \begin{equation}\label{e3.21}
     \limsup_{s\to+\infty}\int^{s+1}_s\int_{|y|\leq\alpha(\tau)}|\nabla w|_{\widetilde{g}}^2dV_{\widetilde{g}}d\tau<\varepsilon_2,
  \end{equation}

or Case2: $\alpha(s)$ satisfies (H2) and
   \begin{equation}\label{e3.22}
     \inf_{s\geq s_3}\int^{s+1}_s\int_{|y|\leq\alpha(\tau)}|\nabla w|_{\widetilde{g}}^2dV_{\widetilde{g}}d\tau<\varepsilon_2.
   \end{equation}
\end{prop}

\noindent\textbf{Proof.} It's suffice to deduce from Proposition \ref{p2.2} that
  $$
   \liminf_{s\to+\infty}a(s+1)\int^{s+1}_s\int_{M}\rho|y|^2|\nabla w|_{\widetilde{g}}^2\varphi^2dV_{\widetilde{g}}d\tau=0
  $$
or
  $$
   \lim_{s\to+\infty}\int^{s+1}_s\int_{M}\rho|y|^2|\nabla w|_{\widetilde{g}}^2\varphi^2dV_{\widetilde{g}}d\tau=0.
  $$
So, conclusion follows from Lemma \ref{l3.4} and Theorem \ref{t3.2}. $\Box$\\

A first one application of Proposition \ref{p3.1} is the following corollary:

\begin{coro}\label{c3.1}
  Let $(M,g)$ be a compact surface equipped with metric $g$ and $u$ be a maximal solution of \eqref{e1.2} on $M\times[0,\omega)$. For any ${p_0}\in M$, let $\varepsilon_2$ and $s_3$ be given in Proposition \ref{p3.1} and set
  $r_0\equiv e^{-s_3/2}$. We have $({p_0},\omega)$ is not a singular point as long as
     \begin{equation}\label{e3.23}
       r^{-2}\int^{\omega-\frac{r^2}{e}}_{\omega-r^2}\int_{D_r({p_0})}|\nabla u|^2dV_gdt<\varepsilon_2
     \end{equation}
  holds for some $0<r\leq r_0$.
\end{coro}

A consequence of Corollary \ref{c3.1} and Vitalli's covering theorem is the following partial regularity result:

\begin{coro}\label{c3.2}
     Letting $(M,g)$ be a compact surface equipped with metric $g$ and $u$ be a maximal solution of \eqref{e1.2} on $M\times[0,\omega)$, we define the blowup set of $u$ by
       $$
        {\mathcal{S}}\equiv\Bigg\{{p_0}\in M\Big|\ r^{-2}\int^{\omega-\frac{r^2}{e}}_{\omega-r^2}\int_{D_r({p_0})}|\nabla u|_g^2dV_gdt\geq\varepsilon_2\ \ \mbox{ holds for all } 0<r<r_0\Bigg\}.
       $$
     Then ${\mathcal{S}}$ contains at most finitely many points.\\
\end{coro}

 Next, we prove the following result of non-degeneracy for finite time blowup:

 \begin{theo}\label{t3.3}
   Under assumptions of Proposition \ref{p3.1}, we take a function $a(s)$ satisfying \eqref{e1.6} and set
      $$
        b(t)=a^{-1}(s),\ \ s=-\log(\omega-t).
      $$
   Then for any ${p_0}\in\Omega$, $({p_0},\omega)$ is not a blowup point for $u$ provided
      \begin{equation}\label{e3.24}
        \limsup_{t\to\omega^-}\sqrt{b(t)(\omega-t)}\sup_{|p-{p_0}|\leq \sqrt{b(t)(\omega-t)}}|\nabla u|_g(p,t)=0.
      \end{equation}
 \end{theo}

\noindent\textbf{Proof.} Under self-similar variables and setting $\alpha(s)=a^{-\frac{1}{2}}(s)$, it yields from \eqref{e3.24}  that
  \begin{equation}\label{e3.25}
     \limsup_{s\to+\infty}\alpha(s)\sup_{|y|\leq \alpha(s)}|\nabla w|_{\widetilde{g}}(y,s)=0.
  \end{equation}
Therefore,
   \begin{equation}\label{e3.26}
     \int^{s+1}_s\int_{|y|\leq\alpha(s)}|\nabla w|_{\widetilde{g}}^2dV_{\widetilde{g}}ds\leq \pi\sup_{\tau\in[s,s+1]}\alpha^2(s)\sup_{|y|\leq\alpha(\tau)}|\nabla w|_{\widetilde{g}}^2(y,\tau)\to0
   \end{equation}
as $s$ tends to infinity. So, $({p_0},\omega)$ is not a singular point by applying Proposition \ref{p3.1} with \eqref{e3.21}. $\Box$\\

As a corollary of Theorem \ref{t3.3}, we have the following result of refined type II singularity:

\begin{coro}\label{c3.3} (Refined type II blowup)
  Under the assumptions of Theorem \ref{t3.3} and taking any function $a(\cdot)$ satisfying \eqref{e1.6}, we set
     $$
      b(t)=a^{-1}(|\log(\omega-t)|).
     $$
  Then for any blow-up point ${p_0}\in M$, we have
    \begin{equation}\label{e3.27}
      \limsup_{t\to\omega^-}\sqrt{b(t)(\omega-t)}\sup_{|p-{p_0}|\leq\sqrt{b(t)(\omega-t)}}|\nabla u|_g(p,t)=+\infty.
    \end{equation}
\end{coro}

\noindent\textbf{Proof.} For any positive non-decreasing function $a(\cdot)$ fulfilling \eqref{e1.6}, we can take another positive non-decreasing function $\overline{a}(\cdot)$ satisfying \eqref{e1.6} and
  $$
    \lim_{s\to+\infty}\frac{\overline{a}(s)}{a(s)}=+\infty.
  $$
Noting that by Theorem \ref{t3.3},
   $$
    \limsup_{t\to\omega^-}\sqrt{\overline{b}(t)(\omega-t)}\sup_{|p-{p_0}|\leq\sqrt{\overline{b}(t)(\omega-t)}}|\nabla u|_g(p,t)>0
   $$
holds for $\overline{b}(t)=\overline{a}^{-1}(|\log(\omega-t)|)$. So, we can conclude that \eqref{e3.27} is also true due to
   $$
    \lim_{t\to\omega^-}\frac{b(t)}{\overline{b}(t)}=+\infty.
   $$
$\Box$

\vspace{40pt}

\section{Finite time blowup when $3\leq m<7$}

In this section, we assume that $M=B_R\subset{\mathbb{R}}^m$ and $N=S^m\subset{\mathbb{R}}^{m+1}$. Given any rotationally symmetric initial data
 \begin{equation}\label{e4.1}
   u_0(x)=\Bigg(\frac{x}{|x|}\sin \theta_0(r),\cos \theta_0(r)\Bigg)
 \end{equation}
with
   $$
     x=(x_1,x_2,\cdots, x_m)\in B_R,
   $$
it's not hard to see that \eqref{e1.2} admits a unique rotationally symmetric solution
  \begin{equation}\label{e4.2}
   u(x,t)=(u_1,u_2,\cdots, u_{m+1})=\Bigg(\frac{x}{|x|}\sin \theta(r,t),\cos \theta(r,t)\Bigg),
  \end{equation}
where $\theta$ satisfies that
   \begin{equation}\label{e4.3}
     \begin{cases}
       \theta_t=\theta_{rr}+\frac{m-1}{r}\theta_r-\frac{m-1}{r^2}\sin \theta\cos \theta, & 0<r<R, t>0\\
       \theta(0,t)=0, \theta(R,t)=b>0, & t>0\\
       \theta(r,0)=\theta_0(r).
     \end{cases}
   \end{equation}

We will prove the following result of finite time blowup and long time existence:

\begin{theo}\label{t4.1}
  Suppose that $3\leq m<7$. We have

  (1) If $b>\vartheta_m$, where $\vartheta_m\in(\frac{\pi}{2},\pi)$ is given by Lemma \ref{p5.1}, then all solutions of \eqref{e4.3} blow up in finite time.

  (2) If $0<b<\frac{\pi}{2}$, then all solutions of \eqref{e4.3} exist for all time.
\end{theo}

Before proving Theorem \ref{t4.1}, we need the following lemma:\\

\begin{lemm}\label{l4.1}
  Let $\theta$ be a maximal solution to \eqref{e4.3} blowing up at finite time $\omega$. Then there exists $0<\omega'<\omega$, such that $\theta_r(0,t)$ is a strictly monotone function in $t\in(\omega',\omega)$.
\end{lemm}

\noindent\textbf{Proof.} Differentiating \eqref{e4.3} in $t$, we know that $v=\theta_t(r,t)$ satisfies \eqref{e5.3} with
   $$
     b(r,t)=(m-1)\cos \theta(r,t),\ \ \ t\in(0,\omega), r\in(0,R].
   $$
So
   $$
     b(0,t)=m-1,\ \ \ b_r(0,t)=0, \ \ \ \forall t\in(0,\omega).
   $$
Furthermore,
   $$
     v(0,t)=v(R,t)\equiv0, \ \ \forall t\in(0,\omega).
   $$
Whenever $v_r(0,t_0)=\theta_{rt}(0,t_0)=0$ at some $t_0\in(0,\omega)$, we have the zero number of $v$ drops at least one. So, there exists some $\omega'\in(0,\omega)$ such that
   $$
     \theta_{rt}(0,t)=v_r(0,t)\not=0, \ \ \forall t\in(\omega',\omega)
   $$
by intersection comparison lemma \ref{l5.1}. Conclusion is drawn. $\Box$\\

Changing $\theta$ to $-\theta$ if necessary, we may assume that
  \begin{equation}\label{e4.4}
   \begin{cases}
     \theta_r(0,t)>0, & \\
     \theta_{rt}(0,t)>0 \ \mbox{ or }\ \theta_{rt}(0,t)<0, & \forall t\in(\omega',\omega).
   \end{cases}
  \end{equation}

\vspace{20pt}

Next, we prove first that the result of finite time blowup:\\

\noindent\textbf{Finite time blowup when $b>\vartheta_m$:}

\noindent{\it Claim 1:} blows up at infinity or finite time: If not, then for a subsequence of time $t_k\to+\infty$,
   \begin{equation}\label{e4.5}
      \theta_k(r,t)\equiv\theta(r,t+t_k)\to \theta_\infty(r,t),\ \ \mbox{ uniformly on } C^{2+\beta,1+\frac{\beta}{2}}([0,R]\times[0,+\infty))
   \end{equation}
as $k\to\infty$, where $\beta\in(0,1)$ and $\theta_\infty$ is also a solution to \eqref{e4.3}. By \eqref{e4.4},
   $$
    \partial_r\theta_\infty(0,t)\equiv constant, \ \ \forall t\geq0
   $$
and thus
   $$
    \partial_{rt}\theta_\infty(0,t)\equiv0, \ \ \forall t\geq0
   $$
hold. So, by Lemma \ref{p5.1}, we derive a steady state
    $$
     \theta_\infty(r,t)\equiv\Phi_\alpha(r),\ \ \forall r\in[0,R], t\geq0
    $$
 of \eqref{e5.6} for some $\alpha>0$. Furthermore, it satisfies that
   $$
    \Phi_\alpha(R)=b>\vartheta_m,
   $$
which contradicts with Lemma \ref{p5.1} below. So, the solution must blow up at infinity or finite time.\\

\noindent{\it Claim 2:} blows up in finite time: If not, then the solution blows up at $t=+\infty$. We need the following lemma for convergence of rescaled profile:

\begin{lemm}\label{l4.2}
  Suppose that the maximal solution of \eqref{e4.3} blows up at infinity, then there exist sequences $\overline{t}_l\to\infty$ and $\overline{\lambda}_l\to0^+$ as $l\to\infty$, such that
    $$
      \lim_{l\to\infty}\theta(\overline{\lambda}_lr, \overline{\lambda}_l^2t+\overline{t}_l)\to\Phi_1(r) \ \mbox{ uniformly in } C^{2+\alpha}_{loc}\Big([0,+\infty)\Big),
    $$
  where $\Phi_a(r)$ is given in Lemma \ref{p5.1} with
     $$
       \partial_r\Phi_a(0)=a.
     $$
\end{lemm}

\noindent\textbf{Proof.} For any $k\in{\mathbb{N}}$, we define $t_k>0$ to be a positive time such that
    \begin{equation}\label{e5.6}
       m(t_k)\sqrt{k-t_k}=\max_{t\in[0,k]}m(t)\sqrt{k-t},
    \end{equation}
where
   $$
    m(t)\equiv\sup_{r\in[0,R]}|\partial_r\theta(r,t)|.
   $$
 Then $\theta$ blows up at infinite time implies that
   \begin{equation}\label{e4.7}
      t_k\to+\infty, \ \ \ m(t_k)\to+\infty\ \ \ \mbox{ as } k\to+\infty,
   \end{equation}
and
    \begin{equation}\label{e4.8}
      \lim_{k\to\infty}m(t_k)\sqrt{k-t_k}=+\infty.
    \end{equation}
Now, rescaling $\theta$ by
   $$
     \theta_k(r,t)\equiv\theta\Big(m^{-1}(t_k)r, m^{-2}(t_k)t+t_k\Big),
   $$
we have
   \begin{eqnarray}\label{e4.9}\nonumber
     \big|\partial_r\theta_k\big|(r,t)&=& m^{-1}(t_k)\big|\partial_r\theta\big|(m^{-1}\big(t_k)r, m^{-2}(t_k)t+t_k\big)\\
     &\leq&\frac{\sqrt{k-t_k}}{\sqrt{k-t_k-m^{-2}(t_k)t}}\to 1 \ \ \ \mbox{ as } k\to\infty,
   \end{eqnarray}
for all $r\in[0,Rm(t_k)], t\in\Big[-t_km^2(t_k),(k-t_k)m^2(t_k)\Big]$. Furthermore, noting that $\theta$ is bounded from above by comparing with some trivial solution $\overline{\theta}(r,t)=\kappa\pi$ with $\kappa\in{\mathbb{N}}$ large, we get
   \begin{eqnarray}\label{e4.10}\nonumber
      \big|\partial_r\theta_k\big|(r,t)&=&m^{-1}(t_k)\big|\partial_r\theta\big|(m^{-1}\big(t_k)r, m^{-2}(t_k)t+t_k\big)\\ \nonumber
     &\leq& C_0 m^{-1}(t_k)\Big(1+(R-m^{-1}(t_k)r)^{-1}+(m(t_k)^{-1}r)^{-1}\Big)\\
     &\leq& C_0\Big(m^{-1}(t_k)+(Rm(t_k)-r)^{-1}+r^{-1}\Big)\to C_0r^{-1}
    \end{eqnarray}
for all $r\in[0,Rm(t_k)], t\in\Big[-t_km^2(t_k),(k-t_k)m^2(t_k)\Big]$ using Lemma \ref{l6.2}.

Consequently, there exists a limiting function $\theta_\infty(r,t)$ satisfying \eqref{e4.3} on $[0,+\infty)\times(-\infty,+\infty)$, such that
   $$
     \theta_k\to\theta_\infty \ \ \mbox{ uniformly on } C^{2+\alpha,1+\frac{\alpha}{2}}_{loc}\Big([0,+\infty)\times(-\infty,+\infty)\Big)
   $$
as $k\to\infty$. As well, there exists a positive constant $M_0$ such that
   \begin{equation}\label{e4.11}
    \begin{cases}
      \sup_{[0,+\infty)\times(-\infty,+\infty)}|\partial_r\theta_\infty|\leq1\\[5pt]
      \sup_{[0,M_0]}|\partial_r\theta_\infty|(\cdot,0)=\sup_{[0,+\infty)}|\partial_r\theta_\infty|(\cdot,0)=1
    \end{cases}
   \end{equation}
by \eqref{e4.10}. Without loss of generality, we also have
   \begin{equation}\label{e4.12}
     \begin{cases}
        \partial_r\theta_\infty(0,t)\geq0,\\[5pt]
        \partial_{rt}\theta_\infty(0,t)\geq0 \mbox{ or } \partial_{rt}\theta_\infty(0,t)\leq0
     \end{cases}
   \end{equation}
for all $t\in{\mathbb{R}}$ by \eqref{e4.4}. Since $\theta_\infty$ is not identical to zero by \eqref{e4.11}, we conclude from \eqref{e4.12} and Corollary \ref{c5.1} that
   \begin{equation}\label{e4.13}
     \begin{cases}
        \partial_r\theta_\infty(0,t)>0,\\[5pt]
        \partial_{rt}\theta_\infty(0,t)>0 \mbox{ or } \partial_{rt}\theta_\infty(0,t)<0
     \end{cases}
   \end{equation}
after some time $t\geq t_0$. In case
    $$
      \partial_{rt}\theta_\infty(0,t)>0\ \ \ \forall t\geq t_0,
    $$
we have
   \begin{equation}\label{e4.14}
      \lim_{t\to+\infty}\partial_r\theta_\infty(0,t)=\alpha\in(0,1].
   \end{equation}
So, one can shift $\theta_\infty$ by
   $$
     \psi_k(r,t)\equiv\theta_\infty(r, t+k)
   $$
for all large $k\in{\mathbb{N}}$, and obtain that
   $$
     \psi_k(r,t)\to\psi_\infty(r,t) \ \ \mbox{ uniformly on } C^{2+\alpha,1+\frac{\alpha}{2}}_{loc}\Big([0,+\infty)\times(-\infty,+\infty)\Big)
   $$
with a limiting function $\psi_\infty$ satisfying \eqref{e4.3} on $[0,+\infty)\times(-\infty,+\infty)$ and
   \begin{equation}\label{e4.15}
       \partial_t\psi_\infty(0,t)=\partial_{rt}\psi_\infty(0,t)\equiv0\ \ \forall t\in{\mathbb{R}}.
   \end{equation}
Now, it follows from Corollary \ref{c5.1} and \eqref{e4.15} that
   \begin{equation}\label{e4.16}
     \psi_\infty(r,t)=\Phi_\alpha(r)\ \ \forall r\in[0,+\infty), t\in{\mathbb{R}}.
   \end{equation}
for some $\alpha>0$. Given any $l\in{\mathbb{N}}$ and taking $k_l\to+\infty$ large enough, after setting
   \begin{eqnarray*}
     \overline{t}_l&\equiv& lm^{-2}(t_{k_l})+t_{k_l},\\
     \overline{\lambda}_l&=& \alpha^{-1}m^{-1}(t_{k_l})
   \end{eqnarray*}
we derive the conclusion of Lemma \ref{l4.2}.

 In case
    $$
      \partial_{rt}\theta_\infty(0,t)<0\ \ \forall t\in{\mathbb{R}},
    $$
 we have
    $$
     \lim_{t\to-\infty}\partial_r\theta_\infty(0,t)=\alpha\in(0,1]
    $$
 and shift $\theta_\infty$ by
    $$
      \psi_k(r,t)\equiv\theta_\infty(r,t-k)
    $$
 for all $k\in{\mathbb{N}}$. A same conclusion follows as above. $\Box$\\

Now, we can complete the proof of finite time blowup by Lemma \ref{l4.2}, \ref{l5.1} and \ref{p5.1} as following:

 By Lemma \ref{l5.1}, there exists a large integer $K$ such that
    \begin{equation}\label{e4.17}
      {\mathcal{Z}}(\theta(\cdot,\overline{t}_l)-\Phi_*)\leq K \ \ \mbox{ for } l=1,2,\cdots
    \end{equation}
 holds for intersection number of $\theta(\cdot,\overline{t}_t)$ with $\Phi_*$ defined in Section 5, where $\{\overline{t}_l\}$ is a sequence of time coming from Lemma \ref{l4.2}. Noting that
    \begin{eqnarray*}
      {\mathcal{Z}}(\theta(\cdot,\overline{t}_l)-\Phi_*)&=&{\mathcal{Z}}\Bigg(\theta(\overline{\lambda}_lr,\overline{t}_l)-\Phi_*(\overline{\lambda}_lr)\Bigg)\\
       &=&{\mathcal{Z}}\Bigg(\theta(\overline{\lambda}_lr,\overline{t}_l)-\Phi_*(r)\Bigg)
    \end{eqnarray*}
by scaling invariant of the singular solution $\Phi_*$ , and
    \begin{equation}\label{e4.18}
      \liminf_{l\to\infty}{\mathcal{Z}}\Big(\theta(\overline{\lambda}_lr,\overline{t}_l)-\Phi_*(r)\Big)\geq{\mathcal{Z}}\Big(\Phi_1-\Phi_*),
    \end{equation}
 we arrive at
   $$
    {\mathcal{Z}}(\Phi_1-\Phi_*)\leq K
   $$
 by Lemma \ref{l4.2}. This contradicts with Lemma \ref{p5.1} and gives the proof of finite time blowup. $\Box$\\

\noindent\textbf{Long time existence when $b<\frac{\pi}{2}$:} Suppose on the contrary, there must be a positive constant $\omega<+\infty$, such that the solution blows up at $\omega$. Thus, there exists a sequence $t_k\to\omega^-$, such that
   \begin{equation}\label{e4.19}
      m(t_k)=\max_{t\in[0,t_k]}m(t)\to+\infty\ \ \mbox{ as } k\to\infty.
   \end{equation}
Re-scaling $\theta$ by
   $$
     \theta_k(r,t)\equiv\theta(m^{-1}(t_k)r,m^{-2}(t_k)t+t_k),
   $$
we obtain that
   \begin{eqnarray}\nonumber\label{e4.20}
     |\partial_r\theta_k|(r,t)&\leq&C_0 m^{-1}(t_k)\Bigg[m(t_k)r^{-1}+\frac{1}{\sqrt{\omega-t_k-m^{-2}(t_k)t}}\Bigg]\\ \nonumber
     &=& C_0\Bigg[r^{-1}+\frac{1}{\sqrt{m^2(t_k)(\omega-t_k)-t}}\Bigg]\\
     &\to& C_0\Bigg[r^{-1}+\frac{1}{\sqrt{T_1-t}}\Bigg]\ \ \mbox{ as } k\to\infty
   \end{eqnarray}
by Lemma \ref{l6.1}, where
   $$
    T_1\equiv\lim_{k\to+\infty}(\omega-t_k)m^2(t_k)\in[0,+\infty].
   $$
Consequently, we get a nontrivial limiting solution $\theta_\infty$ satisfying \eqref{e4.3} on
  $$
   [0,+\infty)\times(-\infty,T_2), \ \ T_2\equiv T_1-C_0^2,
  $$
such that
   $$
    \theta_k\to\theta_\infty\ \ \mbox{ unifromly on } C^{2+\alpha,1+\frac{\alpha}{2}}_{loc}\big([0,+\infty)\times(-\infty,T_2)\big)
   $$
for a subsequence $k=k_j$. Furthermore, \eqref{e4.13} holds for all $t<0$. And thus
   $$
     m^{-1}(t_k)\partial_r\theta\Big(0,(T_2-1)m^{-2}(t_k)+t_k\Big)=\partial_r\theta_k\Big(0,T_2-1\Big)\geq\frac{1}{2}\partial_r\theta_\infty\Big(0,T_2-1\Big)\geq\delta>0
   $$
for $k$ large. As a result, we get
   \begin{equation}\label{e4.21}
      \limsup_{t\to\omega^-}\partial_r\theta(0,t)=+\infty.
   \end{equation}
However, since $b<\frac{\pi}{2}$, after comparing $\theta$ with $\Phi_a$ when $a$ large enough, we have
    \begin{equation}\label{e4.22}
        \theta(r,t)\leq\Phi_a(r),\ \ \ \forall r\in[0,1], t\in[0,\omega)
    \end{equation}
and hence
     \begin{equation}\label{e4.23}
         \partial_r\theta(0,t)\leq a \ \ \forall r\in[0,1], t\in[0,\omega),
     \end{equation}
contradicting with \eqref{e4.21}. So, the solution can not blow up in finite time. $\Box$

\vspace{40pt}

\section{Type I rate when $3\leq m<7$}

In this section, we will prove that the blowup rates of solutions in Theorem \ref{t5.1} are always type I. In fact, we have the following result.

\begin{theo}\label{t5.1}
  Let $\theta$ be the maximal solution of \eqref{e4.3} on $[0,R]\times[0,\omega)$ with $0<\omega<+\infty$. If $3\leq m<7$ we have
     \begin{equation}\label{e5.1}
         \limsup_{t\to\omega^-}(\omega-t)\sup_{0<r\leq1}\Bigg(\frac{m-1}{r^2}\sin^2\theta+\theta_r^2\Bigg)(r,t)<+\infty.
     \end{equation}
\end{theo}

\vspace{10pt}

It's not hard to see that \eqref{e5.1} is equivalent to
  \begin{equation}\label{e5.2}
     \limsup_{t\to\omega^-}(\omega-t)\sup_{0<r\leq R}\theta_r^2(r,t)<+\infty
  \end{equation}
by differential intermediate value theorem and boundary condition $\theta(0,t)=0, \forall t\in[0,\omega)$.

 To prove the above theorem, we need several crucial lemmas under below. The first one is the Sturm-Liouville type theorem for zero number of solution from
\begin{equation}\label{e5.3}
    \frac{\partial v}{\partial t}=v_{rr}+\frac{m-1}{r}v_r-\frac{b(r,t)}{r^2}v, \ 0<r<R, t\in(t_1,t_2),
\end{equation}
where $b(r,t)$ is a bounded function satisfying
  $$
    b(0,t)\equiv m-1, \ b_r(0,t)\equiv 0,\ \ \forall t\in(t_1,t_2)
  $$
in case $v(0,t)\equiv0, \forall t\in(t_1,t_2)$. Furthermore, the following boundary condition
  \begin{equation}\label{e5.4}
    v(r,t)\equiv0 \mbox{ or } v(r,t)\not=0,\ \ \ \forall r=0,R,\ \  t\in(t_1,t_2)
  \end{equation}
is imposed.\\

\begin{lemm}\label{l5.1}
  Let $v$ be a classical solution of \eqref{e5.3} on $[0,R]\times(t_1,t_2)$ which is not identical to zero and satisfies \eqref{e5.4} for some $0<R<+\infty$. We define
    $$
     {\mathcal{Z}}(v(\cdot,t))\equiv\sharp\Big\{r\in[0,R]\Big|\ v(r,t)=0\Big\}
    $$
  to be the zero number of $v(\cdot,t)$ counting the multiplicity. Then

(i) ${\mathcal{Z}}(v(\cdot,t))<\infty$ for any $t_1<t<t_2$,

(ii) ${\mathcal{Z}}(v(\cdot,t))$ is a monotone non-increasing function in time $t$,

(iii) if $v(r_0,t_0)=v_r(r_0,t_0)=0$ for some $0\leq r_0\leq R$ and $t_1<t_0<t_2$, then
    $$
     {\mathcal{Z}}(v(\cdot,t))>{\mathcal{Z}}(v(\cdot,s))\ \ \mbox{ for any } t_1<t<t_0<s<t_2.
    $$
\end{lemm}

The original version of Lemma \ref{l5.1} can be found in \cite{CP} ({see also \cite{A2}) for semilinear heat equation of Fujita type. It was later generalized to the current version for harmonic heat flow in \cite{Du}. It's also notable to remark that when the end point $r=R$ is replaced by a moving free boundary $r=R(t)$, conclusion of Theorem \ref{l5.1} still holds true. One need only using the transformation
   $$
     \overline{v}(r,t)=v(R^{-1}(t)r,t).
   $$
As a consequence of the theorem, we also have the following corollary:

\begin{coro}\label{c5.1}
  Let $v$ be a classical solution of \eqref{e5.3} on $[0,R]\times(t_1,t_2)$ or on $[0,+\infty)\times(t_1,t_2)$, which satisfies \eqref{e5.4} when $0<R<+\infty$. Suppose that for some $t_1<t_*<t^*<t_2$ and $r^*\in[0,R]$ (or $r^*\in[0,+\infty)$ respectively), there holds
     \begin{equation}\label{e5.5}
       v_r(r^*,t)=v(r^*,t)=0\ \ \forall t\in[t_*,t^*],
     \end{equation}
  then $v(r,t)\equiv0$.
\end{coro}

\noindent\textbf{Proof.} If $v$ satisfies \eqref{e5.4} for $0<R<+\infty$, conclusion follows from Lemma \ref{l5.1} since when $v$ not identical to zero, ${\mathcal{Z}}(v(\cdot,t))$ can drop only finitely many zeros and hence contradict with \eqref{e5.5}. In case $v$ is a solution of \eqref{e5.3} on $[0,+\infty)\times(t_1,t_2)$, a same reason can be applied to exclude the possibility of $|v|(r,t)>0$
when $(r,t)$ lies near some $(r_0,t_0)\in(r^*,+\infty)\times(t_*,t^*)$.
In fact, if not, then $v$ must be identical to zero in $[0,r_0)\times(t_*,t^*)$ by Lemma \ref{l5.1}. This contradicts with our assumption $|v|(r,t)>0$ near $(r_0,t_0)$. Therefore, $v(r,t)\equiv0$ for $r\geq0, t\in (t_*,t^*)$. The proof was done. $\Box$\\

Next, let's consider the half-entire solution of
  \begin{equation}\label{e5.6}
  \begin{cases}
    \tau(\Phi)\equiv \Phi_{rr}+\frac{m-1}{r}\Phi_r-\frac{m-1}{r^2}\sin \Phi\cos \Phi=0, &  r\in[0,+\infty)\\
    \Phi(0)=0, \Phi(r)>0, & \forall r>0
  \end{cases}
  \end{equation}
on half line. Setting $\Phi_1$ to be a solution of \eqref{e5.6} with
   $$
     \Phi'_1(0)=1,
   $$
we have $\Phi_a(r)=\Phi_1(ar)$ is also a solution of \eqref{e5.6} with
   $$
     \Phi'_a(0)=a.
   $$
It's notable also that the trivial singular solution to \eqref{e5.6} is given by
   $$
     \Phi_*(r)\equiv\frac{\pi}{2}.
   $$
We have the following lemma concerning the number of intersection points between $\Phi_a$ and $\Phi_*$:

\begin{prop}\label{p5.1}
  Suppose that $3\leq m<7$, we have
     $$
      {\mathcal{Z}}(\Phi_a-\Phi_*)=+\infty
     $$
  and
     $$
       \lim_{r\to+\infty}\Phi_a(r)=\frac{\pi}{2},\ \ \ \ \lim_{r\to+\infty}|\Phi'_a(r)|=0
     $$
  for any $a>0$. Furthermore,
     $$
       \max_{r\in[0,+\infty)}\Phi_a(r)=\vartheta_m\in\Big(\frac{\pi}{2},\pi\Big)
     $$
  for any $a>0$.
\end{prop}

It is remarkable Biernat-Seki have shown in \cite{BS} the result in Lemma \ref{p5.1} does not hold for $m\geq7$ since $\Phi_a(r)$ is monotone increasing for all $r>0$.\\

We divide the proof into three lemmas:

\begin{lemm}\label{l5.2}
  For any $k\in{\mathbb{N}}$, there exist a decreasing sequence $\omega_{2k-1}\in(\frac{\pi}{2},\pi)$ and a increasing sequence $\omega_{2k}\in(0,\frac{\pi}{2})$ such that $\Phi$ increases strictly from $0$ to $\omega_1$, then decreases strictly from $\omega_1$ to $\omega_2$, and so on.
\end{lemm}

\noindent\textbf{Proof.} Letting $\Phi$ be a solution to \eqref{e5.6} satisfying $\Phi'(0)>0$, there exists a maximal interval $[0,\omega_1)$ with $\omega\leq+\infty$, such that $\Phi$ increases strictly until $\Phi=\omega_1$. In case $\omega_1<+\infty$, we have
  \begin{equation}\label{e5.7}
     \liminf_{\Phi\to\omega_1^-}\Phi_{rr}\leq0, \ \ \lim_{\Phi\to\omega_1^-}\Phi_{r}=0.
  \end{equation}
Regarding $r$ as a function of $\Phi$, one gets
   $$
    \Phi_r=\frac{1}{r'},\ \ \Phi_{rr}=-\frac{1}{(r')^3}r''.
   $$
Setting $y=y(\Phi)\equiv\frac{r'}{r}$, equation \eqref{e5.6} changes to
   \begin{equation}\label{e5.8}
      \begin{cases}
        y'-(m-2)y^2+\frac{m-1}{2}y^3\sin 2\Phi=0, & \forall\Phi>0,\\
        y(0)=+\infty, y(\omega_1)=+\infty, 0<y(\Phi)<+\infty, & \forall 0<\Phi<\omega_1.
      \end{cases}
   \end{equation}
Letting $w_1=y^{-2}$ for $\Phi\in(0,\omega_1)$, there holds
   \begin{equation}\label{e5.9}
     \begin{cases}
       w_1'+2(m-2)\sqrt{w_1}-(m-1)\sin2\Phi=0, & \forall\Phi\in(0,\omega_1)\\
       w_1(0)=0, w_1(\omega_1)=0, 0<w_1(\Phi)<+\infty, & \forall \Phi\in(0,\omega_1).
     \end{cases}
   \end{equation}
Setting $\Phi_0\equiv\Phi(r_0)$ for $r_0\in(0,r_1)$, where $r_1\in(0,\infty]$ is the first time $\Phi$ reaching $\omega_1$, we have the resolution formula
   \begin{equation}\label{e5.10}
     \log r-\log r_0=\int^{\Phi(r)}_{\Phi_0}\frac{d\tau}{\sqrt{w_1(\tau)}}
   \end{equation}
for each $r\in(0,r_1)$.

\noindent\textbf{\it Claim 1: when $\omega_1<+\infty$ and $\omega_1\not=\frac{k\pi}{2}, k\in{\mathbb{N}}$, there holds
      \begin{equation}\label{e5.11}
         r_1<+\infty.
      \end{equation}
In case $\omega_1=\frac{k\pi}{2}$ for some $k\in{\mathbb{N}}$, we have
     \begin{equation}\label{e5.12}
       r_1=+\infty.
     \end{equation}
}

\noindent\textbf{Proof of Claim 1.} Actually, when $\omega_1<+\infty$ and $\omega_1\not=\frac{k\pi}{2}, k\in{\mathbb{N}}$, the vanishing order of $w_1$ near $\omega_1$ must be one by equation \eqref{e5.9}. Thus, it follows from the resolution formula \eqref{e5.10} that $r_1<+\infty$. Another hand, when $\omega_1=\frac{k\pi}{2}$ for some $k\in{\mathbb{N}}$, the vanishing order of $w_1$ near $\omega_1$ must be two by equation \eqref{e5.9}. As a result, it follows from \eqref{e5.10} that $r_1=+\infty$. $\Box$\\

\noindent\textbf{\it Claim 2: $\omega_1\geq\pi/2$.}

\noindent\textbf{Proof of Claim 2.} If not, then $\omega_1<\pi/2$. By Claim 1, there must be $r_1<+\infty$. However, it follows from \eqref{e5.6} and \eqref{e5.7} that $\sin(2\omega_1)\leq0$. Contradiction holds. $\Box$\\

\noindent\textbf{\it Claim 3: for $3\leq m<7$, there holds $\omega_1>\frac{\pi}{2}$.}

 \noindent\textbf{Proof of Claim 3.} By Claim 1, if $\omega_1=\frac{\pi}{2}$, the solution $\Phi(r)$ increases for all $r>0$ and tends to $\omega_1=\pi/2$ as $r\to\infty$.

Noting that $\overline{w}_1\equiv\sin^2(2\Phi)$ is a solution to
   $$
    \begin{cases}
      \overline{w}_1'+2(m-2)\sqrt{\overline{w}_1}-(m-1)\sin2\Phi\\
      \ \ \ \ \ \ \ \ \ \ \ =\sin2\Phi\Big[4\cos2\Phi+(m-3)\Big], & \forall \Phi\in(0,\pi/2)\\
      \overline{w}_1(0)=\overline{w}_1(\pi/2)=0
    \end{cases}
   $$
for $m<7$, which implies that $\cos2\Phi+\frac{m-3}{4}$ changes sign exactly once, we get
    \begin{equation}\label{e5.13}
      w_1(\Phi)\leq\sin^2(2\Phi), \ \ \forall \Phi\in(0,\pi/2).
    \end{equation}
Substituting into \eqref{e5.9} yields that
    \begin{equation}\label{e5.14}
      w_1(\Phi)\leq(m-3)\cos^2\Phi, \ \ \forall\Phi\in(0,\pi/2).
    \end{equation}
Setting
  $$
    F(\Phi)\equiv-2(m-2)\sqrt{w_1}+(m-1)\sin(2\Phi), \ \ F(\pi/2)=0,
  $$
we want to show that $F$ is positive near $\pi/2$ on left hand side. Since $F(0)$, this is equivalent to prove that
   \begin{eqnarray*}
    F'(\Phi)&-&-(m-2)\frac{w'_1}{\sqrt{w_1}}+2(m-1)\cos(2\Phi)\\
     &=&-(m-2)\Bigg[-2(m-2)+(m-1)\frac{\sin(2\Phi)}{\sqrt{w_1}}\Bigg]+2(m-1)\cos(2\Phi)\\
     &\leq&-(m-2)\Bigg[-2(m-2)+(m-1)\frac{\sin(2\Phi)}{\sqrt{m-3}\cos\Phi}\Bigg]+2(m-1)\cos(2\Phi)\\
     &\sim&-(m-2)\Bigg[-2(m-2)+\frac{2(m-1)}{\sqrt{m-3}}\Bigg]-2(m-1)<0, \ \ \forall \Phi\in(\pi/2-\delta,\pi/2)
   \end{eqnarray*}
for $m<7$, where \eqref{e5.14} has been used. So, $w_1$ must be increasing near $\Phi=\pi/2$, which contradicts with the fact $w_1(\Phi)>0, \forall\Phi\in(0,\pi/2)$ and $w_1(\pi/2)=0$. The proof of Claim 3 was done. $\Box$\\

\noindent\textbf{\it Claim 4: $\omega_1<\pi$.}

\noindent\textbf{Proof.} Noting that $v=\sin^2\Phi$ satisfies that
   \begin{equation}\label{e5.15}
     \begin{cases}
       v'+2(m-2)\sqrt{v}-(m-1)\sin2\Phi=2(m-2)\sin\Phi(1-\cos\Phi)\geq0, & \forall \Phi\in[0,\pi]\\
       v(0)=0, v(\pi)=0, 0<v(\Phi)<+\infty, & \forall \Phi\in(0,\pi),
     \end{cases}
   \end{equation}
we have
   $$
     w_1(\Phi)<\sin^2\Phi\ \ \ \forall \Phi\in(0,\min\{\omega_1,\pi\})
   $$
by comparing $v$ with $w_1$. So, $\omega_1<\pi$. $\Box$\\

 Since $\Phi'=0$ when $\Phi$ reaches $\omega_1$, by \eqref{e5.6}, we know that $\Phi$ will decrease strictly below $\omega_1$ due to $\omega_1\in(\frac{\pi}{2},\pi)$. Now, defining
 $\omega_2<\omega_1$ such that
    \begin{equation}\label{e5.16}
      \begin{cases}
         y'-(m-2)y^2+\frac{m-1}{2}y^3\sin2\Phi=0, & \Phi\in(\omega_2,\omega_1),\\
         y(\omega_2)=-\infty, y(\omega_1)=-\infty, -\infty<y(\Phi)<0, & \forall \Phi\in(\omega_2,\omega_1)
      \end{cases}
    \end{equation}
 with $y=\frac{r'}{r}<0$, we have $\omega_2<\frac{\pi}{2}$ by \eqref{e5.6}. Setting $w_2=y^{-2}$ as above, one gets similarly that
    \begin{equation}\label{e5.17}
      \begin{cases}
         w_2'-2(m-2)\sqrt{w_2}-(m-1)\sin2\Phi=0, & \forall \Phi\in(\omega_2,\omega_1)\\
         w_2(\omega_2)=0, w_2(\omega_1)=0, 0<w_2(\Phi)<+\infty, & \forall\Phi\in(\omega_2,\omega_1),
      \end{cases}
    \end{equation}
 which is different from \eqref{e5.9} with negative second term. Now, comparing $w_2$ with $w_1$ on $[\omega_2,\omega_1]$, we get $\omega_2\in(0,\frac{\pi}{2})$.

   Similarly, we can define $\omega_3>\frac{\pi}{2}>\omega_2$ such that $w_3\equiv\frac{r^2}{(r')^2}$ satisfies \eqref{e5.9} on $(\omega_2,\omega_3)$. After comparing $w_3$ with $w_1$, we have $\omega_3<\omega_1$. So, Lemma \ref{l5.2} holds true by a bootstrap argument.

\begin{lemm}\label{l5.3}
  The half-entire solution must increase and decrease infinitely many times when $m\geq3$.
\end{lemm}

 In fact, for any $k\in{\mathbb{N}}$,
   \begin{equation}\label{e5.18}
    \Big(\log r\Big)'=\begin{cases}
        \frac{1}{\sqrt{w_{2k-1}(\Phi)}}, & \Phi\in(\omega_{2k-2},\omega_{2k-1})\\
        -\frac{1}{\sqrt{w_{2k}(\Phi)}}, & \Phi\in(\omega_{2k},\omega_{2k-1}).
     \end{cases}
   \end{equation}
 Therefore, given $K\in{\mathbb{N}}$, we have
    \begin{equation}\label{e5.19}
      \begin{cases}
        \log\Big(r(\omega_{2K})\Big)=\lim_{\Phi_0\to0^+}\Bigg(\log\Big(\frac{\Phi_0}{\Phi'(0)}\Big)+\int^{\omega_1}_{\Phi_0}\frac{d\Phi}{\sqrt{w_1(\Phi)}}\Bigg)+\Bigg(\Sigma_{k=2}^K\int^{\omega_{2k-1}}_{\omega_{2k}}\frac{d\Phi}{\sqrt{w_{2k}(\Phi)}}+\Sigma_{k=1}^K\int^{\omega_{2k-1}}_{\omega_{2k}}\frac{d\Phi}{\sqrt{w_{2k-1}(\Phi)}}\Bigg)\\[5pt]
        \log\Big(r(\omega_{2K+1})\Big)=\lim_{\Phi_0\to0^+}\Bigg(\log\Big(\frac{\Phi_0}{\Phi'(0)}\Big)+\int^{\omega_1}_{\Phi_0}\frac{d\Phi}{\sqrt{w_1(\Phi)}}\Bigg)+\Bigg(\Sigma_{k=2}^{K+1}\int^{\omega_{2k+1}}_{\omega_{2k}}\frac{d\Phi}{\sqrt{w_{2k}(\Phi)}}+\Sigma_{k=1}^K\int^{\omega_{2k-1}}_{\omega_{2k}}\frac{d\Phi}{\sqrt{w_{2k-1}(\Phi)}}\Bigg)
      \end{cases}
    \end{equation}
 Another hand, for any $k$, by \eqref{e5.9} for $w_{2k-1}$ or \eqref{e5.17} for $w_{2k}$,
     \begin{equation}\label{e5.20}
      \begin{cases}
        w'_{2k-1}(r)\not=0 & \mbox{ for } r=\omega_{2k-2}, \omega_{2k-1}\\
        w'_{2k}(r)\not=0 & \mbox{ for } r=\omega_{2k-1}, \omega_{2k}.
      \end{cases}
     \end{equation}
 A combination of \eqref{e5.19} and \eqref{e5.20} yields
     \begin{equation}\label{e5.21}
       r(\omega_{K})<+\infty
     \end{equation}
 for any $K\geq1$. Lemma \ref{l5.3} holds true. (It's notable that this claim is not true for $m=2$)

\begin{lemm}\label{l5.4}
   \begin{equation}\label{e5.22}
     \lim_{K\to+\infty}r(\omega_K)=+\infty.
   \end{equation}
\end{lemm}

If not, then there exists $0<r_{max}<+\infty$ such that
   $$
     r(\omega_K)\uparrow r_{max} \ \ \mbox{ as } K\to+\infty.
   $$
Noting that
   $$
     \Phi'(r(\omega_K))=0 \ \ \forall K,
   $$
there holds
   $$
     \Phi(r_{max})=\frac{\pi}{2}, \Phi'(r_{max})=\Phi''(r_{max})=0.
   $$
Then $\Phi$ must be identical to $\frac{\pi}{2}$, contradicting with $\Phi(0)=0$. So, Lemma \ref{l5.4} holds true.

To complete the proof of lemma, note first that it follows from Lemma \ref{l5.2} that
   $$
    0\leq w_k(\Phi)\leq \sin^2\Phi\leq1, \ \ \mbox{ for }\begin{cases}
       \Phi\in(\omega_{2K-2},\omega_{2K-1}), &  k=2K-1,\\
       \Phi\in(\omega_{2K},\omega_{2K-1}), & k=2K.
    \end{cases}
   $$
Therefore,
   $$
     |\Phi'(r)|=\frac{\sqrt{w_k(r)}}{r}\leq\frac{1}{r}\to0
   $$
as $r\to+\infty$. Finally, we show that
   \begin{equation}\label{e5.23}
     \lim_{k\to\infty}\omega_{2k}=\lim_{k\to\infty}\omega_{2k+1}=\frac{\pi}{2}.
   \end{equation}
Suppose on the contrary, since $\omega_{2k}$ is monotone increasing and lies below $\frac{\pi}{2}$, $\omega_{2k+1}$ is monotone decreasing and lies above $\frac{\pi}{2}$, we have
   $$
     \lim_{k\to\infty}\omega_{2k}=\omega_*<\omega^*=\lim_{k\to\infty}\omega_{2k+1}.
   $$
Since $w_{k}$ lies in a bounded set of $C^1([\omega_*,\omega^*])\cap C^2((\omega_*,\omega^*))$, passing to the limits, we get two limiting functions $w^{+}_\infty$ and $w^{-}_\infty$ satisfying
  \begin{equation}\label{e5.24}
     (w^{(\pm)}_{\infty})'\pm2(m-2)\sqrt{w^{(\pm)}_{\infty}}-(m-1)\sin2\Phi=0, \ \ \  \forall \Phi\in(\omega_*,\omega^*)
  \end{equation}
and
  \begin{equation}\label{e5.25}
    w^{\pm}_\infty(\omega_*)=w^{\pm}_\infty(\omega^*)=0.
  \end{equation}
Integrating over $(\omega_*,\omega^*)$, we get
   $$
     2(m-1)\int^{\omega^*}_{\omega_*}\sqrt{w^{\pm}_\infty}d\Phi=\pm(m-1)\int^{\omega^*}_{\omega_*}2\sin2\Phi d\Phi,
   $$
which is impossible since L.H.S. is positive. So \eqref{e5.23} holds true, and all conclusions of the lemma is drawn. $\Box$\\

Let $\theta$ be a maximal solution to \eqref{e5.3} on $[0,R)\times[0,\omega)$. For any $0<t<\omega$, we define
   \begin{equation}\label{e5.26}
     m(t)\equiv\sup_{r\in[0,R]}|\theta_r|(r,t).
   \end{equation}
Then we have the following characterization for profile of type II blowup which will be proven in Section 6:

\begin{prop}\label{p5.2}
   Let $\theta$ be a maximal solution to \eqref{e4.3} on $[0,R)\times[0,\omega)$. Suppose that the blowup at $\omega<+\infty$ is type II, then there exist two sequences $t_l\to\omega^-$ and $\lambda_l\to0^+$ for $l=1,2,\cdots$, such that
      \begin{equation}\label{e5.27}
        \theta\Big(\lambda_lr,t_l\Big)\to\Phi_1(r) \ \mbox{ or } -\Phi_1(r),\ \ \ \mbox{ as } l\to+\infty
      \end{equation}
   uniformly on any compact set of $[0,+\infty)$.\\
\end{prop}

Now, let's turn to prove our main theorem \ref{t5.1}. Suppose that it is not true, then the rate of blowup is type II. Without loss of generality, we may assume that the limiting function is given by $\Phi_1(r)$ in \eqref{e5.27}. By Lemma \ref{l5.1}, there exists a large integer $K$ such that
  \begin{equation}\label{e5.28}
    {\mathcal{Z}}(\theta(\cdot,t_l)-\Phi_*)\leq K \ \ \mbox{ for } l=1,2,\cdots,
  \end{equation}
where $\{t_l\}_{l=1}^\infty$ is a sequence of time coming from Proposition \ref{p5.2}. Noting that
  \begin{eqnarray*}
    {\mathcal{Z}}(\theta(\cdot,t_l)-\Phi_*)&=&{\mathcal{Z}}\Bigg(\theta(\lambda_lr,t_l)-\Phi_*(\lambda_lr)\Bigg)\\
      &=&{\mathcal{Z}}\Bigg(\theta(\lambda_lr,t_l)-\Phi_*(r)\Bigg)
  \end{eqnarray*}
since singular solution $\Phi_*$ is invariant under scaling, and
   $$
    \liminf_{l\to+\infty}{\mathcal{Z}}\Bigg(\theta(\lambda_lr,t_l)-\Phi_*(r)\Bigg)\geq{\mathcal{Z}}(\Phi_1-\Phi_*),
   $$
we conclude that from \eqref{e5.28} that
   $$
     {\mathcal{Z}}(\Phi_1-\Phi_*)\leq K,
   $$
which contradicts with Lemma \ref{p5.1}. So, Theorem \ref{t5.1} holds true. $\Box$

\vspace{40pt}

\section{Singular profiles of type II blowup}

In this section, we will verify the validity of Proposition \ref{p5.2}. Let's start with a-priori estimate on derivative of the solution to \eqref{e5.3}:

\begin{lemm}\label{l6.1}
  Let $\theta$ be a solution to \eqref{e4.3} on $[0,R]\times[0,T), R\geq1$. There exists a constant $C_0>0$ depending only on $m, ||\theta_{0r}||_{L^\infty[0,R]}$ and $T$, such that
    \begin{equation}\label{e6.1}
       |\theta_r(r,t)|\leq C_0\Big(r^{-1}+(T-t)^{-\frac{1}{2}}\Big), \ \ \forall r\in(0,R], 0\leq t<T.
    \end{equation}
\end{lemm}

\noindent\textbf{Proof.} Setting
   $$
     u(x,t)=\Bigg(\frac{x}{|x|}\sin\theta(|x|,t),\cos\theta(|x|,t)\Bigg),
   $$
we have $u$ satisfies that
   \begin{equation}\label{e6.2}
     \begin{cases}
       u_t=\triangle u+|\nabla u|^2u, & (x,t)\in B_R\times(0,\omega),\\
       u(x,t)=\Big(\frac{x}{R}\sin b,\cos b\Big),&  x\in\partial B_R, t\in[0,\omega),\\
       u(x,0)=u_0(x)\equiv\Big(\frac{x}{|x|}\sin\theta_0(|x|),\cos\theta_0(|x|)\Big).
     \end{cases}
   \end{equation}
For any $a\in\overline{B_R}$, rescaling $u$ by
  $$
    u(x,t)=w\Bigg(\frac{x-a}{\sqrt{T-t}},-\log(T-t)\Bigg),
  $$
it's clear that $w=w_{(a,T)}(y,s)$ satisfies
   \begin{eqnarray}\label{e6.3}
    &w_s-\triangle w+\frac{1}{2}y\cdot\nabla w=|\nabla w|^2w,&\\ \nonumber &\forall y\in\Omega_s\equiv\{y\in{\mathbb{R}}^m|\ \ |a+e^{-\frac{s}{2}}y|\leq R\}, s\geq-\log T.&
   \end{eqnarray}
Multiplying \eqref{e6.3} by $\rho w_s, \rho(y)\equiv e^{-\frac{|y|^2}{4}}$ and performing integration by parts, one gets
  \begin{eqnarray}\label{e6.4}\nonumber
    \frac{1}{2}\frac{d}{ds}\int_{\Omega_s}\rho|\nabla w|^2&=&-\int_{\Omega_s}\rho w_s^2-\frac{1}{4}e^{-s}\int_{\partial\Omega_s}\rho(y\cdot\nu)\theta_r^2d\sigma_y\\ \nonumber
    &&+\frac{m-1}{4R^2}e^{-s}\int_{\partial\Omega_s}\rho(y\cdot\nu)\sin^2\theta d\sigma_y\\
    &\leq&-\int_{\Omega_s}\rho w_s^2-\frac{1}{4}e^{-s}\int_{\partial\Omega_s}\rho(y\cdot\nu)\theta_r^2d\sigma_y+Ce^{-s}.
  \end{eqnarray}
So, after integrating over time, there holds
   \begin{equation}\label{e6.5}
       \int_{B_1\cap\Omega_s}|\nabla w_{(a,T)}|^2dy\leq C_1
   \end{equation}
for all $a\in\overline{B_R}, s\geq-\log T$ and some positive constant $C_1$ depending only on $||\theta_{0r}||_{L^\infty([0,R])}$ and $T$ (Under below, we will also denote $C_2, C_3$ etc to be positive constants depending only on $||\theta_{0r}||_{L^\infty([0,R])}$ and $T$). Using the relation
   $$
     w_{(a,T)}(y,s)=w_{(0,T)}(y+e^{\frac{s}{2}}a,s),
   $$
it's inferred from \eqref{e6.5} that
   \begin{equation}\label{e6.6}
      \int_{y\in\Omega_s,\ |y-e^{\frac{s}{2}}a|\leq1}|\nabla w_{(0,T)}|^2dy\leq C_1, \ \ \forall a\in\overline{B_R}.
   \end{equation}
Noting that
   \begin{eqnarray*}
     w_{(0,T)}(y,s)&=& u(e^{-\frac{s}{2}}y,T-e^{-s})\\
       &=&\Bigg(\frac{y}{|y|}\sin\theta(e^{-\frac{s}{2}}|y|,T-e^{-s}), \cos\theta(e^{-\frac{s}{2}}|y|,T-e^{-s})\Bigg),
   \end{eqnarray*}
if one sets
   $$
     \Theta(r,s)\equiv\theta(e^{-\frac{s}{2}}r,T-e^{-s})
   $$
we can obtain
   \begin{equation}\label{e6.7}
      \int_{|r|\leq e^{\frac{s}{2}}R, |r-e^{\frac{s}{2}}r_0|\leq3/4}\Bigg(\frac{m-1}{r^2}\sin^2\Theta+\Theta_r^2\Bigg)dr\leq C_1
   \end{equation}
for all $r_0\in[0,R]$ and $s\in[s_0,s_0+1], s_0\geq-\log T$. Because $\Theta(r,s)$ satisfies
  \begin{equation}\label{e6.8}
     \Theta_s=\Theta_{rr}+\frac{m-1}{r}\Theta_r-\frac{1}{2}r\Theta_r-\frac{m-1}{2r^2}\sin2\Theta,
  \end{equation}
after setting
   $$
    \overline{\Theta}(\varrho,\tau)=\Theta(r,s),\ \ \ \varrho=re^{-\frac{s-s_0}{2}}, \tau=-e^{-(s-s_0)},
   $$
we have
   \begin{equation}\label{e6.9}
     \overline{\Theta}_\tau=\overline{\Theta}_{\varrho\varrho}+\frac{m-1}{\varrho}\overline{\Theta}_{\varrho}-\frac{m-1}{2\varrho^2}\sin2\overline{\Theta}
   \end{equation}
and
   \begin{equation}\label{e6.10}
      \int_{|\varrho-\varrho_0|\leq1/2}\overline{\Theta}_{\varrho}^2d\varrho\leq 2C_1, \ \forall\varrho_0\geq1, -1\leq\tau\leq-e^{-1}.
   \end{equation}
Using the parabolic estimate for one dimensional equation \eqref{e6.9} under \eqref{e6.10}, it yields that
   \begin{equation}\label{e6.11}
      \Big|\overline{\Theta}_{\varrho}\Big|(\varrho,\tau)\leq C_2, \ \forall\varrho\geq1, -1\leq\tau\leq -e^{-1}.
   \end{equation}
Writing back to $\theta$, one gets
   \begin{equation}\label{e6.12}
      \Big|\Theta_r\Big|(r,s)\leq C_3, \ \forall r\geq1, s\geq-\log T.
   \end{equation}

 To complete the proof of the lemma, we need only to estimate $w=w_{(0,T)}(y,s)$ near origin. In fact, taking any $0<T'<T$ and noting that
    \begin{eqnarray*}
      \theta(r,t)&=& \Theta_{(0,T)}\Bigg(\frac{r}{\sqrt{T-t}},-\log(T-t)\Bigg)\\
        &=& \Theta_{(0,T')}\Bigg(\frac{r}{\sqrt{T'-t}},-\log(T'-t)\Bigg),
    \end{eqnarray*}
 we have
     $$
       \Theta_{(0,T)}(r,s)=\Theta_{(0,T')}\Bigg(r\sqrt{\frac{e^{-s}}{T'-T+e^{-s}}}, s+\log\frac{e^{-s}}{T'-T+e^{-s}}\Bigg).
     $$
Recalling \eqref{e6.12} for $\Theta=\Theta_{(0,T')}$ and setting
  $$
    \lambda=\frac{e^{-s}}{T'-T+e^{-s}},
  $$
 we conclude that
   \begin{eqnarray}\label{e6.13}\nonumber
     \Big|\frac{\partial}{\partial r}\Theta_{(0,T)}\Big|(r,s)&=&\sqrt{\lambda}\Big|\frac{\partial}{\partial r}\Theta_{(0,T')}\Big|\Bigg(r\sqrt{\lambda}, s+\log\lambda\Bigg)\\
     &\leq& C_3(1+r^{-1}),\ \ \forall r>0, s\geq-\log T
   \end{eqnarray}
by taking
    $$
      T'=T-e^{-s}+r^2e^{-s}.
    $$
Therefore,
   \begin{eqnarray*}
     \Big|\frac{\partial\theta}{\partial r}\Big|(r,t)&=&\frac{1}{\sqrt{T-t}}\Big|\frac{\partial\Theta_{(0,T)}}{\partial r}\Big|\Bigg(\frac{r}{\sqrt{T-t}},-\log(T-t)\Bigg)\\
     &\leq& C_3(r^{-1}+(T-t)^{-\frac{1}{2}}),
   \end{eqnarray*}
and the proofs were done. $\Box$\\

For eternal solution of \eqref{e4.3}, we have the following result of a-priori bound:\\

\begin{lemm}\label{l6.2}
  Let $\theta$ be a bounded solution to \eqref{e4.3} on $[0,R]\times[0,+\infty)$ for some $R>0$. There exists a constant $C=C(m, \theta, R)>0$, such that
    \begin{equation}\label{e6.14}
      |\theta_r(r,t)|\leq C\Big(1+r^{-1}+(R-r)^{-1}\Big), \ \forall r\in [0,R], t\geq0.
    \end{equation}
\end{lemm}

\noindent\textbf{Proof.} Suppose on the contrary, then for any $k\in{\mathbb{N}}$, there exists a sequence of $(x_k,t_k)\in[0,R)\times[0,k)$, such that
    \begin{eqnarray}\nonumber\label{e6.15}
       &&\frac{|\theta_r(r_k,t_k)|}{1+r_k^{-1}+(R-r_k)^{-1}+(k-t_k)^{-1/2}}\\
       && \ \ \ \ \ \ \ \ =\sup_{(r,t)\in[0,R)\times[0,k)}\frac{|\theta_r(r,t)|}{1+r^{-1}+(R-r)^{-1}+(k-t)^{-1/2}}\to+\infty
    \end{eqnarray}
as $k$ large. It's clear that
   \begin{eqnarray}\nonumber\label{e6.16}
     t_k&\to&\infty, \ m_k\equiv|\theta_r(r_k,t_k)|\to\infty,\\
      r_km_k&\to&\infty,\ (R-r_k)m_k\to\infty,\ (k-t_k)^{1/2}m_k\to\infty.
   \end{eqnarray}
Now, re-scaling $\theta$ by
   $$
    \theta_k(r,t)\equiv\theta(m_k^{-1}r+r_k, m_k^{-2}t+t_k),
   $$
we have
   \begin{eqnarray}\nonumber\label{e6.17}
     |\partial_r\theta_k(r,t)|&=& m_k^{-1}|\partial_r\theta|(m_k^{-1}r+r_k, m_k^{-2}t+t_k),\\ \nonumber
     &\leq&\frac{1+(m_k^{-1}r+r_k)^{-1}+(R-r_k-m_k^{-1}r)^{-1}+(k-t_k-m_k^{-2}t)^{-1/2}}{1+r_k^{-1}+(R-r_k)^{-1}+(k-t_k)^{-1/2}},\\ \nonumber
     &\leq&\frac{m_k^{-1}+(r+r_km_k)^{-1}+((R-r_k)m_k-r)^{-1}+((k-t_k)m_k^{-2}-t)^{-1/2}}{m_k^{-1}+(r_km_k)^{-1}+((R-r_k)m_k)^{-1}+((k-t_k)m_k)^{-1/2}},\\
     &\to& 1, \ \ \forall x\in[-r_km_k, (R-r_k)m_k], t\in[(-t_km_k^2, k-t_k)m_k^2]
   \end{eqnarray}
and
   \begin{equation}\label{e6.18}
     \partial_t\theta_k-\partial^2_r\theta_k=\frac{m-1}{r+r_km_k}\partial_r\theta_k-\frac{m-1}{(r+r_km_k)^2}\sin\theta_k\cos\theta_k.
   \end{equation}
Therefore, after passing to the limits, one gets a limiting solution $\theta_\infty$ on ${\mathbb{R}}^2$ satisfying that
  \begin{equation}\label{e6.19}
    \begin{cases}
      \partial_t\theta_\infty-\partial^2_r\theta_\infty=0, & (r,t)\in{\mathbb{R}}^2,\\
       \partial_r\theta_{\infty}(0,0)=\pm1, \ |\partial_r\theta_{\infty}|(r,t)\leq1, & \forall (r,t)\in{\mathbb{R}}^2.
    \end{cases}
  \end{equation}
Since $\partial_r\theta_\infty$ is also a solution to heat equation which attains its maximum or minimum at $(0,0)$, it follows from strong maximum principle that $\partial_r\theta_\infty(r,t)\equiv\pm1$ on ${\mathbb{R}}^2$. This fact contradicts with the bounded assumption of $\theta$. Conclusion is drawn. $\Box$\\

Now, we turn to prove Proposition \ref{p5.2} similar as Lemma \ref{l4.2}:\\

\noindent\textbf{Proof of Proposition \ref{p5.2}:} For any $k\in{\mathbb{N}}$, let's define $t_k\in(0,\omega)$ by
   \begin{equation}\label{e6.20}
      m(t_k)\sqrt{\omega-\frac{1}{k}-t_k}=\max_{t\in[0,\omega-\frac{1}{k}]}m(t)\sqrt{\omega-\frac{1}{k}-t}.
   \end{equation}
Since $\theta$ blows up at $\omega$ with type II rate, we have
   \begin{equation}\label{e6.21}
     t_k\to\omega, \ \ m(t_k)\to+\infty\ \ \ \mbox{ as } k\to+\infty,
   \end{equation}
and
    \begin{equation}\label{e6.22}
      \lim_{k\to\infty} m(t_k)\sqrt{\omega-\frac{1}{k}-t_k}=+\infty.
    \end{equation}
Rescaling $\theta$ by
     $$
       \theta_k(r,t)\equiv\theta\Big(m^{-1}(t_k)r, m^{-2}(t_k)t+t_k\Big),
     $$
there holds
     \begin{eqnarray}\label{e6.23}\nonumber
       \big|\partial_r\theta_k\big|(r,t)&=& m^{-1}(t_k)\big|\partial_r\theta\big|(m^{-1}(t_k)r, m^{-2}(t_k)t+t_k)\\
       &\leq&\frac{\sqrt{\omega-\frac{1}{k}-t_k}}{\sqrt{\omega-\frac{1}{k}-t_k-m^{-2}(t_k)t}}\to1\ \ \ \mbox{ as } k\to\infty
     \end{eqnarray}
for all $r\in[0,Rm(t_k)], t\in[-t_km^2(t_k), \omega-\frac{1}{k}-t_k-m^{-2}(t_k)t]$.

We still also have
   \begin{eqnarray}\nonumber\label{e6.24}
     |\partial_r\theta_k|(r,t)&\leq&C_0 m^{-1}(t_k)\Bigg[m(t_k)r^{-1}+\frac{1}{\sqrt{\omega-t_k-m^{-2}(t_k)t}}\Bigg]\\ \nonumber
     &=& C_0\Bigg[r^{-1}+\frac{1}{\sqrt{m^2(t_k)(\omega-t_k)-t}}\Bigg]\\
     &\to& C_0r^{-1}\ \ \mbox{ as } k\to\infty
   \end{eqnarray}
by Lemma \ref{l6.1}. Consequently, there exists a limiting function $\theta_\infty(r,t)$ satisfying \eqref{e5.3} on $[0,+\infty)\times(-\infty,+\infty)$, such that
     $$
       \theta_k\to\theta_\infty \ \mbox{ uniformly on } C^{2+\alpha,1+\frac{\alpha}{2}}_{loc}\big([0,+\infty)\times(-\infty,+\infty)\big)
     $$
 as $k\to\infty$. Furthermore, there exists a positive constant $M_0$ such that
    \begin{equation}\label{e6.25}
      \begin{cases}
         \sup_{[0,+\infty)\times(-\infty,+\infty)}|\partial_r\theta_\infty|\leq1,\\
         \sup_{[0,M_0]}|\partial_r\theta_\infty|(\cdot,0)=\sup_{[0,+\infty)}|\partial_r\theta_\infty|(\cdot,0)=1
      \end{cases}
    \end{equation}
 by \eqref{e6.24}. Without loss of generality, we assume that
   \begin{equation}\label{e6.26}
      \begin{cases}
        \partial_r\theta_\infty(0,t)>0,\\
        \partial_{rt}\theta_\infty(0,t)>0 \mbox{ or } \partial_{rt}\theta_\infty(0,t)<0
      \end{cases}
   \end{equation}
 by Lemma \ref{l4.1}. Now, a same argument as in Lemma \ref{l4.2} shows that Proposition \ref{p5.2} holds true. $\Box$

\vspace{40pt}

\section{Harmonic heat flow from $S^m$ when degree is no less than $2$}

When considering harmonic heat flow from $S^m$ to $S^m\subset{\mathbb{R}}^{m+1}$, the metric $g$ under stereo polar coordinates of sphere is given by
   $$
     g_{ij}(x)=\frac{1}{(1+|x|^2)^2}\delta_{ij}.
   $$
Therefore, when $u_0(x)\in S^m\subset{\mathbb{R}}^{m+1}$ is symmetric under rotation and inversion, it's not hard to see that the solution $u(x,t)$ is also symmetric under rotation and inversion, and satisfies that
   \begin{equation}\label{e7.1}
     \begin{cases}
       \frac{\partial u}{\partial t}=g^{-1}(|x|)(\triangle u+|\nabla u|^2u),& (x,t)\in B_1\times(0,\omega),\\
        u(x,0)=u_0(x),\\
        u(x,t)\Big|_{\partial B_1}=u_0\Big|_{\partial B_1}
     \end{cases}
   \end{equation}
by uniqueness, where $g(|x|)=\frac{1}{(1+|x|^2)^2}$. So, when $u_0(x)$ takes the form \eqref{e4.1}, we have $u(x,t)$ takes the form of \eqref{e4.2}, and the equation of $u$ is reduced to
   \begin{equation}\label{e7.2}
     \begin{cases}
       \theta_t=g^{-1}(r)\Big(\theta_{rr}+\frac{m-1}{r}\theta_r-\frac{m-1}{2r^2}\sin2\theta\Big), & r\in(0,1), t\in(0,\omega),\\
       \theta(r,0)=\theta_0(r),\\
       \theta(0,t)=0, \theta(1,t)=b.
     \end{cases}
   \end{equation}
Also, if the degree of $u_0$ is no less than two, we have $b\geq\pi>\vartheta_m$.

Noting that $\frac{1}{4}\leq g(r)\leq 1, r\in[0,1]$, we will prove the following result on finite time type I blowup when dimension greater than two:

\begin{theo}\label{t7.1}
  Consider the rotatory-inversion symmetric harmonic heat flow from $S^m$ to $S^m\subset{\mathbb{R}}^{m+1}, 3\leq m<7$. If $u_0$ is a super-harmonic map whose degree is no less than $2$, we have the solution blows up in finite time with type I rate.\\
\end{theo}

To complete the proof of Theorem \ref{t7.1}, we need several lemmas. The first one is the existence of super-harmonic map under our settings.

\begin{lemm}\label{l7.1}
  For any $b\geq\pi$, there exist infinitely many solutions $\theta\equiv\theta_0(r)$ satisfying
    \begin{equation}\label{e7.3}
      \begin{cases}
        \theta_{rr}+\frac{m-1}{r}\theta_r-\frac{m-1}{2r^2}\sin2\theta\geq0, & \forall r\in(0,1)\\
        \theta(0)=0, \theta(1)=b.
      \end{cases}
    \end{equation}
\end{lemm}

\noindent\textbf{Proof.} By Lemma \ref{p5.1}, there exist infinitely many $a>0$ such that
  \begin{equation}\label{e7.4}
    \Phi_a(1)\in\Big(\frac{\pi}{2},\pi\Big)
  \end{equation}
in case of $3\leq m<7$. If one choose $\theta(r)=\Phi_a(r)+\beta r$ for $\beta=b-\Phi_a(1)>0$, we have
  \begin{eqnarray*}
    z&\equiv&\theta_{rr}+\frac{m-1}{r}\theta_r-\frac{m-1}{2r^2}\sin2\theta\\
    &=&\partial^2_r\Phi_a+\frac{m-1}{r}\Phi_a+\frac{m-1}{r}\beta-\frac{m-1}{2r^2}\sin(2\Phi_a+2\beta r)\\
    &=&\frac{m-1}{r}\beta-\frac{m-1}{r^2}\cos(2\Phi_a+\beta r)\sin(\beta r).
  \end{eqnarray*}
Thus, it's clear that $z>0$ when
  $$
   \beta r\in[2k\pi, (2k+1)\pi], \ \ 2\Phi_a+\beta r\in\Bigg[\Bigg(2l+\frac{1}{2}\Bigg)\pi, \Bigg(2l+\frac{3}{2}\Bigg)\pi\Bigg]
  $$
or
  $$
   \beta r\in[(2k+1)\pi, (2k+2)\pi], \ \ 2\Phi_a+\beta r\in\Bigg[2l\pi,\Bigg(2l+\frac{1}{2}\Bigg)\pi\Bigg]\bigcup\Bigg[\Bigg(2l+\frac{3}{2}\Bigg)\pi, 2(l+1)\pi\Bigg],
  $$
for nonnegative integers $k,l$. In case of
  $$
   \beta r\in[2k\pi, (2k+1)\pi], \ \ 2\Phi_a+\beta r\in\Bigg[2l\pi,\Bigg(2l+\frac{1}{2}\Bigg)\pi\Bigg]\bigcup\Bigg[\Bigg(2l+\frac{3}{2}\Bigg)\pi, 2(l+1)\pi\Bigg],
  $$
we have
  \begin{eqnarray*}
   z&\geq&\frac{m-1}{r}\beta-\frac{m-1}{r^2}\sin(\beta r-2k\pi)\\
    &\geq&\frac{m-1}{r}\beta-\frac{m-1}{r^2}(\beta r-2k\pi)\geq0
  \end{eqnarray*}
Similarly, for
  $$
   \beta r\in[(2k+1)\pi, (2k+2)\pi], \ \ 2\Phi_a+\beta r\in\Bigg[\Bigg(2l+\frac{1}{2}\Bigg)\pi, \Bigg(2l+\frac{3}{2}\Bigg)\pi\Bigg],
  $$
we still have
  \begin{eqnarray*}
   z&\geq&\frac{m-1}{r}\beta+\frac{m-1}{r^2}\sin(\beta r)\\
    &\geq&\frac{m-1}{r}\beta-\frac{m-1}{r^2}(\beta r-(2k+1)\pi)\geq0.
  \end{eqnarray*}
The proof was done. $\Box$\\

\begin{lemm}\label{l7.2}
  For any initial datum $\theta_0$ from Lemma \ref{l7.1}, the solution $\theta(r,t)$ of \eqref{e7.2} is monotone non-decreasing in time. Furthermore, it is a super-solution to
   \begin{equation}\label{e7.5}
     \theta_t=g_*\Bigg(\theta_{rr}+\frac{m-1}{r}\theta_r-\frac{m-1}{2r^2}\sin2\theta\Bigg), \ \ r\in(0,1), t>0
   \end{equation}
 and a sub-solution to
    \begin{equation}\label{e7.6}
     \theta_t=g^*\Bigg(\theta_{rr}+\frac{m-1}{r}\theta_r-\frac{m-1}{2r^2}\sin2\theta\Bigg), \ \ r\in(0,1), t>0
   \end{equation}
 under a same initial-boundary condition, where
   $$
    g_*\equiv\inf_{r\in[0,1]} g^{-1}(r), \ \ g^*\equiv\sup_{r\in[0,1]}g^{-1}(r).
   $$
\end{lemm}

\noindent\textbf{Proof.} Taking derivative of $t$ on \eqref{e7.2}, one deduces that the equation of $v\equiv\theta_t$ is given by
  $$
   \begin{cases}
     v_t=g^{-1}(r)\Big(v_{rr}+\frac{m-1}{r}v_r-\frac{m-1}{r^2}v\cos2\theta\Big) , & r\in(0,1), t>0,\\
     v(r,0)\geq0, &  \forall r\in[0,1]\\
     v(0,t)=v(1,t)=0, & \forall t\geq0.
   \end{cases}
  $$
By maximum principle, we conclude that $v(r,t)$ is nonnegative everywhere. Thus, $\theta$ is a super-solution to \eqref{e6.5} and a sub-solution to \eqref{e6.6}. $\Box$\\

\noindent\textbf{Proof of finite time blowup:} When
 $$
   u_0(x)=\Bigg(\frac{x}{|x|}\sin\theta_0(r), \cos\theta_0(r)\Bigg)
 $$
is a super-harmonic map with degree no less than two, it's not hard to see that $\theta_0$ satisfies \eqref{e6.3} with some $b\geq\pi$. Transforming the solution $\theta_*$ of \eqref{e6.5} by $\theta_*(r,g_*t)$ and using Theorem \ref{t4.1} for \eqref{e4.3}, we conclude that $\theta_*$ blows up in finite time $\omega_*\in(0,+\infty)$. Therefore, there is a time $\omega\in(0,\omega_*]$ such that the solution $\theta$ of \eqref{e7.2} blows up at $\omega$ by maximum principle. This complete the first part of Theorem \ref{t7.1} for finite time blowup. $\Box$\\

  To show the second part of type I blowup rate, let's continue our lemmas. The third one is a parallel version of Lemma \ref{l7.1}.

\begin{lemm}\label{l7.3}
  Let $\theta$ be a solution to \eqref{e7.2} on $[0,1]\times[0,T)$. Under assumption of Theorem \ref{t7.1}, there exists a constant $C_0>0$ depending only on $m, ||\theta_{0r}||_{L^\infty[0,1]}$ and $T$, such that
    \begin{equation}\label{e7.7}
       |\theta_r(r,t)|\leq C_0\Big(r^{-1}+(T-t)^{-\frac{1}{2}}\Big), \ \ \forall r\in(0,R], 0\leq t<T.
    \end{equation}
\end{lemm}

\noindent\textbf{Proof.} Noting that under assumption of Theorem \ref{t7.1}, $\widetilde{\theta}(r,t)\equiv\theta(r,g^*t)$ satisfies that
  \begin{equation}\label{e7.8}
   \widetilde{\theta}_t\leq \widetilde{\theta}_{rr}+\frac{m+1}{r}\widetilde{\theta}_r-\frac{m-1}{2r^2}\sin2\widetilde{\theta}, \ \ r\in(0,1), t>0.
  \end{equation}
Minor change in proof of Lemma \ref{l6.1} shows that \eqref{e6.1} for solution of \eqref{e4.3} also holds for solution of \eqref{e6.8}. Written back to $\theta$ yields the desire inequality \eqref{e6.7}. $\Box$\\

The next proposition is a parallel version of Proposition \ref{p5.2}.

\begin{prop}\label{p7.1}
   Let $\theta$ be a maximal solution to \eqref{e7.2} on $[0,1)\times[0,\omega)$. Suppose that the blowup at $\omega<+\infty$ is type II, then there exist two sequences $t_l\to\omega^-$ and $\lambda_l\to0^+$ for $l=1,2,\cdots$, such that
      \begin{equation}\label{e7.9}
        \theta\Big(\lambda_lr,t_l\Big)\to\Phi_1(r) \ \mbox{ or } -\Phi_1(r),\ \ \ \mbox{ as } l\to+\infty
      \end{equation}
   uniformly on any compact set of $[0,+\infty)$.\\
\end{prop}

\noindent\textbf{Proof.} The proof is very similar to proof of Proposition \ref{p5.2}, expect replacing Lemma \ref{l6.1} by Lemma \ref{l7.3} and replacing Lemma \ref{l4.1} (in deduction of \eqref{e6.26}) by Lemma \ref{l7.2} to deduce that
   \begin{equation}\label{e7.10}
    \begin{cases}
     \theta_r(0,t)>0, & \forall t\geq0\\
     \theta_{rt}(0,t)\geq0, & \forall t\geq0,
    \end{cases}
   \end{equation}
where we have used
  $$
   \theta_t(0,t)\equiv0, \ \ \theta_t(r,t)\geq0, \ \ \forall r\in[0,1], t\geq0.
  $$
The conclusion of the proposition was drawn. $\Box$\\

The main technical difficulty in proof of Theorem \ref{t7.1} is no suitable zero comparison lemma on hand so far. However, when considering the equation
  \begin{equation}\label{e7.11}
   \frac{\partial v}{\partial t}=g^{-1}(r)\Bigg(v_{rr}+\frac{m-1}{r}v_r-\frac{b(r,t)}{r^2}v\Bigg), \ \ 0<r<1, t\in(t_1,t_2)
  \end{equation}
under the boundary condition
  \begin{equation}\label{e7.12}
    v(0,t)\not=0, \ \ v(1,t)\equiv0 \mbox{ or } v(1,t)\not=0,\ \ \forall t\in(t_1,t_2),
  \end{equation}
we still have the following zero comparison lemma by applying a result of Matano  (Theorem 1 in \cite{Ma}, see also  Proposition 2 in \cite{A1}).

\begin{lemm}\label{l8.4}
  Consider the solution to \eqref{e6.11} under boundary condition \eqref{e6.12}, we have
   \begin{equation}\label{e7.13}
     {\mathcal{Z}}(v(\cdot,t))\geq{\mathcal{Z}}(v(\cdot,s)), \ \ \forall t_1<t<s<t_2,
   \end{equation}
  where
    $$
     {\mathcal{Z}}(v(\cdot,t))\equiv\sharp\Big\{r\in[0,1]\Big|\ v(r,t)=0\Big\}
    $$
  as before.
\end{lemm}

\noindent\textbf{Proof.} A key observation is that under boundary condition \eqref{e7.12}, the zeros are located away from origin. Therefore, the equation \eqref{e7.11} becomes a linear second order parabolic equation with continuous coefficient. So, one can impose the zero comparison theorem of Matano \cite{Ma} to our case. $\Box$\\

Now, we can complete the proof of Theorem \ref{t7.1}.

\noindent\textbf{Proof of type I rate.} The proof is similar to that of Theorem \ref{t5.1} by comparing $\theta(\lambda_kr, t_l)$ with $\Phi_*$ and using Lemma \ref{p5.1}. $\Box$

\vspace{40pt}

\section*{Acknowledgments}
The author would like to express his deepest gratitude to Professors Xi-Ping Zhu, Kai-Seng Chou, Xu-Jia Wang and Neil Trudinger for their constant encouragements and warm-hearted helps. This paper is also dedicated to the memory of Professor Dong-Gao Deng.

\end{document}